\documentclass[reqno,centertags, 12pt]{amsart}
\usepackage{amsmath,amsthm,amscd,amssymb}
\usepackage{latexsym}
\newcounter{smalllist}
\newenvironment{SL}{\begin{list}{{\rm\roman{smalllist})}}{%
\setlength{\topsep}{0mm}\setlength{\parsep}{0mm}\setlength{\itemsep}{0mm}%
\setlength{\labelwidth}{2em}\setlength{\leftmargin}{2em}\usecounter{smalllist}%
}}{\end{list}}

\newcommand{\bigtimes}{\mathop{\mathchoice%
{\smash{\vcenter{\hbox{\LARGE$\times$}}}\vphantom{\prod}}%
{\smash{\vcenter{\hbox{\Large$\times$}}}\vphantom{\prod}}%
{\times}%
{\times}%
}\displaylimits}

\newcommand{\bbC}{{\mathbb{C}}}
\newcommand{\bbD}{{\mathbb{D}}}
\newcommand{\bbE}{{\mathbb{E}}}

\newcommand{\bbT}{{\mathbb{T}}}
\newcommand{\bbU}{{\mathbb{U}}}
\newcommand{\bbZ}{{\mathbb{Z}}}

\newcommand{\calC}{{\mathcal C}}
\newcommand{\calE}{{\mathcal E}}
\newcommand{\calG}{{\mathcal G}}
\newcommand{\calH}{{\mathcal H}}
\newcommand{\calI}{{\mathcal I}}

\newcommand{\calL}{{\mathcal L}} 
\newcommand{\calM}{{\mathcal M}}
 
\newcommand{\calP}{{\mathcal P}}

\newcommand{\dott}{\,\cdot\,}

\newcommand{\lb}{\label}
\newcommand{\f}{\frac}

\newcommand{\ol}{\overline}
\newcommand{\ti}{\tilde  }

\newcommand{\tr}{\text{\rm{Tr}}}
\newcommand{\dist}{\text{\rm{dist}}}

\newcommand{\ch}{\text{\rm{ch}}}
\newcommand{\intt}{\text{\rm{int}}}

\newcommand{\ess}{\text{\rm{ess}}}

\newcommand{\s}{\text{\rm{s}}}

\newcommand{\supp}{\text{\rm{supp}}}

\newcommand{\bi}{\bibitem}

\newcommand{\beq}{\begin{equation}}
\newcommand{\eeq}{\end{equation}}
\newcommand{\ba}{\begin{align}}
\newcommand{\ea}{\end{align}}
\newcommand{\veps}{\varepsilon}

\DeclareMathOperator{\Real}{Re}
\DeclareMathOperator{\ran}{Ran}
\DeclareMathOperator{\Ima}{Im}

\DeclareMathOperator*{\slim}{s-lim}

\let\det=\undefined\DeclareMathOperator{\det}{det}
\allowdisplaybreaks
\numberwithin{equation}{section}

\newtheorem{theorem}{Theorem}[section]

\newtheorem{lemma}[theorem]{Lemma}
\newtheorem{corollary}[theorem]{Corollary}
\theoremstyle{definition}

\newtheorem{conjecture}[theorem]{Conjecture}

\theoremstyle{remark}

\newcommand{\abs}[1]{\lvert#1\rvert}

\begin{document}
\title[OPUC: New Results]{Orthogonal Polynomials on \\ the Unit Circle: New Results}
\author[B.~Simon]{Barry Simon}

\thanks{$^*$ Mathematics 253-37, California Institute of Technology, Pasadena, CA 91125. 
E-mail: bsimon@caltech.edu. Supported in part by NSF grant DMS-0140592} 

\date{May 5, 2004}

\begin{abstract} We announce numerous new results in the theory of orthogonal 
polynomials on the unit circle. 
\end{abstract}

\maketitle

\section{Introduction} \lb{s1} 

I am completing a comprehensive look at the theory of orthogonal polynomials on 
the unit circle (OPUC; we'll use OPRL for the real-line case). These two  
$500+$-page volumes \cite{OPUC1,OPUC2} to appear in the same AMS series that 
includes Szeg\H{o}'s celebrated 1939 book \cite{Szb} contain numerous new 
results. Our purpose here is to discuss the most significant of these new 
results. Besides what we say here, some joint new results appear instead in 
papers with I.~Nenciu \cite{NenSim}, Totik \cite{SimTot}, and Zlato\v{s} \cite{SZprep}. 
We also note that some of the results I discuss in this article are unpublished 
joint work with L.~Golinskii (Section~3.2) and with Denisov (Section 4.2). Some 
other new results appear in \cite{Haifaproc}. 

Throughout, $d\mu$ will denote a nontrivial (i.e., not a finite combination of 
delta functions) probability measure on $\partial\bbD$, the boundary of $\bbD=
\{z\mid \abs{z}<1\}$. We'll write 
\begin{equation} \lb{1.1} 
d\mu(\theta)=w(\theta)\, \f{d\theta}{2\pi} + d\mu_\s (\theta) 
\end{equation}
where $d\mu_\s$ is singular and $w\in L^1 (\partial\bbD, \f{d\theta}{2\pi})$. 

Given $d\mu$, one forms the monic orthogonal polynomials, $\Phi_n (z;d\mu)$,  
and orthonormal polynomials 
\begin{equation} \lb{1.2} 
\varphi_n (z;d\mu)= \f{\Phi_n(z;d\mu)}{\|\Phi_n\|_{L^2}}
\end{equation}
If one defines 
\begin{equation} \lb{1.3} 
\alpha_n = -\ol{\Phi_{n+1}(0)} 
\end{equation}
then the $\Phi$'s obey a recursion relation 
\begin{equation} \lb{1.4} 
\Phi_{n+1}(z)=z\Phi_n(z) -\bar\alpha_n \Phi_n^*(z) 
\end{equation} 
where $\,{}^*\,$ is defined on degree $n$ polynomials by 
\begin{equation} \lb{1.5} 
P_n^*(z)=z^n \,\ol{P_n (1/\bar z)} 
\end{equation}
\eqref{1.4} is due to Szeg\H{o} \cite{Szb}. The cleanest proofs are in Atkinson 
\cite{Atk64} and Landau \cite{Land}. The $\alpha_n$ are called Verblunsky 
coefficients after \cite{V35}. Since $\Phi_n^*$ is orthogonal to $\Phi_{n+1}$, 
\eqref{1.4} implies 
\begin{align} 
\|\Phi_{n+1}\|^2 &= (1-\abs{\alpha_n}^2) \|\Phi_n \|^2  \lb{1.6} \\
&=\prod_{j=0}^n (1-\abs{\alpha_j}^2) \lb{1.7} 
\end{align} 

It is a fundamental result of Verblunsky \cite{V35} that $\mu\mapsto 
\{\alpha_n\}_{n=0}^\infty$ sets up a one-one correspondence between nontrivial 
probability measures and $\bigtimes_{n=0}^\infty \bbD$. 

A major focus in the book \cite{OPUC1,OPUC2} and in our new results is the view 
of $\{\alpha_n\}_{n=0}^\infty \leftrightarrow \mu$ as a spectral theory problem 
analogous to the association of $V$\! to the spectral measure $-\f{d}{dx^2}+V(x)$ 
or of Jacobi parameters to a measure in the theory of OPRL. 

We divide the new results in major sections: Section~2 involving relations to 
Szeg\H{o}'s theorem, Section~3 to the CMV matrix, Section~4 on miscellaneous 
results, Section~5 on the case of periodic Verblunsky coefficients, and 
Section~6 to some interesting spectral theory results in special classes of 
Verblunsky coefficients. 

\medskip
I'd like to thank P.~Deift, S.~Denisov, L.~Golinskii, S.~Khruschchev, 
R.~Killip, I.~Nenciu, P.~Nevai, F.~Peherstorfer, V.~Totik, and A.~Zlato\v{s} 
for useful discussions.

\bigskip
\section{Szeg\H{o}'s Theorem} \lb{s2} 

In the form first given by Verblunsky \cite{V36}, this says, with $\mu$ given by 
\eqref{1.1}, that   

\begin{equation} \lb{2.1} 
\prod_{j=0}^\infty (1-\abs{\alpha_j}^2)=\exp \biggl(\int_0^{2\pi} \log (w(\theta)) \, 
\f{d\theta}{2\pi}\biggr) 
\end{equation}

\subsection{Szeg\H{o}'s Theorem via Entropy}\lb{s2.1} 

The sum rules of Killip-Simon \cite{KS} can be viewed as an OPRL analog of \eqref{2.1} 
so, not surprisingly, \eqref{2.1} has a ``new" proof that mimics that in \cite{KS}. 
Interestingly enough, while the proof in \cite{KS} has an easy half that depends on 
semicontinuity of the entropy and a hard half (that even after simplifications in 
\cite{SZ,Sim288} is not so short), the analog of the hard half for \eqref{2.1} follows 
in a few lines from Jensen's inequality and goes back to Szeg\H{o} in 1920 \cite{Sz20,Sz21}. 
Here's how this analogous proof goes (see \cite[Section~2.3]{OPUC1} for details):  
\begin{SL} 
\item[(a)] (well-known, goes back to Szeg\H{o} \cite{Sz20,Sz21}). By \eqref{1.7}, 
\begin{align} 
\prod_{j=0}^n (1-\abs{\alpha_j}^2) &\geq \int \exp[ \log (w(\theta)) + 
\log \abs{\Phi_n^*(e^{i\theta})}^2]\, \f{d\theta}{2\pi} \lb{2.2}  \\
&\geq \exp\biggl( \int \log (w(\theta)) + 2\log \abs{\Phi_n^*(e^{i\theta})}\biggr) 
\f{d\theta}{2\pi} \lb{2.3} \\
&= \exp\biggl(\int_0^{2\pi} \log (w(\theta)) \, \f{d\theta}{2\pi}\biggr) \lb{2.4} 
\end{align} 
where \eqref{2.2} uses $d\mu \geq w(\theta)\f{d\theta}{2\pi}$, \eqref{2.3} is Jensen's 
inequality, \and \eqref{2.4} uses the fact that since $\Phi_n^*$ is nonvanishing in 
$\bar\bbD$, $\log \abs{\Phi_n^*(z)}$ is harmonic there and $\Phi_n^*(0)=1$.  

\item[(b)] The map $d\mu\mapsto\int_0^{2\pi} \log (w(\theta))\f{d\theta}{2\pi}$ is a 
relative entropy and so weakly upper semicontinuous in $\mu$ by a Gibbs' variational 
principle: 
\begin{equation} \lb{2.5} 
\int_0^{2\pi} \log(w(\theta))\, \f{d\theta}{2\pi} = 
\inf_{\substack{f\in C(\partial\bbD)\\ f>0}}\, \biggl[ \int f(\theta)\, d\mu(\theta) -1 
-\int \log(f(\theta))\, \f{d\theta}{2\pi}\biggr]
\end{equation} 

\item[(c)] By a theorem of Geronimus \cite{Ger44}, if 
\begin{equation} \lb{2.6}  
d\mu_n(\theta) = \f{d\theta}{2\pi \abs{\varphi_n(e^{i\theta})}^2} 
\end{equation}
(the Bernstein-Szeg\H{o} approximations),  then $d\mu_n\to d\mu$ weakly and the 
Verblunsky coefficients of $d\mu_n$ obey 
\begin{equation} \lb{2.7} 
\alpha_j (d\mu_n) =\begin{cases} \alpha_j(d\mu) & j=0, \dots, n-1 \\
0 & j\geq n \end{cases} 
\end{equation}
Therefore, by the weak upper semicontinuity of (b), 
\begin{equation} \lb{2.8}  
\int_0^{2\pi} \log (w(\theta))\, d\mu \geq\limsup_{n\to\infty}  \int_0^{2\pi} 
-\log (\abs{\varphi_n (e^{i\theta})}^2)\, \f{d\theta}{2\pi}
\end{equation} 

\item[(d)] Since $\abs{\varphi_n (e^{i\theta})}=\abs{\varphi_n^*(e^{i\theta})}= 
\prod_{j=0}^{n-1} (1-\abs{\alpha_j}^2)^{-1/2} \abs{\Phi_n^*(e^{i\theta})}$, 
the same calculation that went from \eqref{2.3} to \eqref{2.4} shows 
\begin{equation} \lb{2.9}  
\exp \biggl[\int_0^{2\pi} -\log (\abs{\varphi_n (e^{i\theta})}^2)\, \f{d\theta}{2\pi} 
\biggr] = \prod_{j=0}^{n-1} (1-\abs{\alpha_j}^2) 
\end{equation}
\eqref{2.4}, \eqref{2.8}, and \eqref{2.9} imply \eqref{2.1} and complete the sketch 
of this proof. 
\end{SL}

\smallskip
We put ``new" in front of this proof because it is closely related to the 
almost-forgotten proof of Verblunsky \cite{V36} who, without realizing he was dealing 
with an entropy or a Gibbs' principle, used a formula close to \eqref{2.5} in his 
initial proof of \eqref{2.1} 

The interesting aspect of this entropy proof is how $d\mu_\s$ is handled en passant --- 
its irrelevance is hidden in \eqref{2.5}.

\subsection{A Higher-Order Szeg\H{o} Theorem}\lb{s2.2} 

\eqref{2.1} implies 
\begin{equation} \lb{2.10}  
\sum_{j=0}^\infty \, \abs{\alpha_j}^2 <\infty \Leftrightarrow \int_0^{2\pi} \log(w(\theta))\, 
\f{d\theta}{2\pi} >-\infty 
\end{equation} 

The following result of the same genre is proven as Theorem~2.8.1 in \cite{OPUC1}: 

\begin{theorem}\lb{T2.1} For any Verblunsky coefficients $\{\alpha_j\}_{j=0}^\infty$, 
\begin{equation} \lb{2.11}  
\sum_{j=0}^\infty\, \abs{\alpha_{j+1}-\alpha_j}^2 + \sum_{j=0}^\infty  
\, \abs{\alpha_j}^4 <\infty \Leftrightarrow \int_0^{2\pi} (1-\cos(\theta)) 
\log(w(\theta))\, \f{d\theta}{2\pi} >-\infty 
\end{equation} 
\end{theorem}  

The proof follows the proof of \eqref{2.10} using the sum rule 
\begin{equation} \lb{2.12}  
\begin{split}
\exp \biggl( -\tfrac12\, \abs{\alpha_0}^2 - \Real (\alpha_0) & + \tfrac12 \sum_{j=0}^\infty \, 
 \abs{\alpha_{j+1}-\alpha_j}^2 \biggr) \prod_{j=0}^\infty (1-\abs{\alpha_j}^2) 
e^{\abs{\alpha_j}^2} \\ 
& = \exp\biggl( \int_0^{2\pi} (1-\cos(\theta)) 
\log(w(\theta))\, \f{d\theta}{2\pi}\biggr)
\end{split}
\end{equation}
in place of \eqref{2.1} The proof of \eqref{2.12} is similar to the proof of \eqref{2.1} 
sketched in Section~\ref{s2.1}. For details, see \cite[Section~2.8]{OPUC1}. 

\smallskip
Earlier than this work, Denisov \cite{Deni} proved that when the left side of \eqref{2.11} 
is finite, then $w(\theta) >0$ for a.e.~$\theta$. In looking for results like \eqref{2.10}, 
we were motivated in part by attempts of Kupin \cite{Kup1,Kup2} and Latpev et al.~\cite{LNS} 
to extend the OPUC results of Killip-Simon (see also \cite{NPVY}). After Theorem~2.1 
appeared in a draft of \cite{OPUC1}, Denisov-Kupin \cite{DenKup2} and Simon-Zlatos 
\cite{SZprep} discussed higher-order analogs.  

\subsection{Relative Szeg\H{o} Function}\lb{s2.3} 

In the approach to sum rules for OPRL called step-by-step, a critical role is played 
by the fact that if $m$ is the $m$-function for a Jacobi matrix, $J$, and $m_1$ is 
the $m$-function for $J_1$, the matrix obtained from $J$ by removing one row and 
column, then 
\begin{equation} \lb{2.13}  
\f{\Ima m_1 (E+ i 0)}{\Ima m(E+i0)}=\abs{a_1 m(E+i0)}^2 
\end{equation}

The most obvious analog of the $m$-function for OPUC is the Carath\'eodory function 
\begin{equation} \lb{2.14}  
F(z) =\int \f{e^{i\theta}+z}{e^{i\theta}-z} \, d\mu(\theta)
\end{equation}
If $\{\alpha_j\}_{j=0}^\infty$ are the Verblunsky coefficients of $d\mu$, the analog of 
$m_1$ is obtained by letting $\beta_j=\alpha_{j+1}$ and $d\mu_1$ the measure with $\alpha_j 
(d\mu_1)=\beta_j$ and $d\mu_1=w_1(\theta)\f{d\theta}{2\pi}+ d\mu_{1,\s}$. 

For $\f{d\theta}{2\pi}$-a.e.~$\theta\in\partial\bbD$, $F(e^{i\theta})\equiv 
\lim_{r\uparrow 1} F(re^{i\theta})$ has a limit and 
\begin{equation} \lb{2.15}  
w(\theta) = \Real F(e^{i\theta}) 
\end{equation}
Thus, as in \eqref{2.13}, we are interested in $\Real F(e^{i\theta})/\Real F_1(e^{i\theta})$ 
which, unlike \eqref{2.13}, is not simply related to $F(e^{i\theta})$. Rather, there is a 
new object $(\delta_0 D)(z)$ which we have found whose boundary values have a magnitude 
equal to the square root of $\Real F(e^{i\theta})/\Real F_1(e^{i\theta})$. 

To define $\delta_0 D$, we recall the Schur function, $f$, of $d\mu$ is defined by 
\begin{equation} \lb{2.16}  
F(z)=\f{1+zf(z)}{1-zf(z)} 
\end{equation}
$f$ maps $\bbD$ to $\bbD$ and \eqref{2.14}/\eqref{2.16} set up a one-one correspondence 
between such $f$'s and probability measures on $\partial\bbD$.  

$\delta_0 D$, the relative Szeg\H{o} function, is defined by 
\begin{equation} \lb{2.17}  
(\delta_0 D)(z)=\f{1-\bar\alpha_0 f(z)}{\rho_0}\, 
\f{1-zf_1}{1-zf}
\end{equation} 
where $f_1$ is the Schur function of $d\mu_1$. One has the following: 

\begin{theorem}\lb{T2.2} Let $d\mu$ be a nontrivial probability measure on $\partial\bbD$ 
and $\delta_0 D$ defined by \eqref{2.17}. Then 
\begin{SL} 
\item[{\rm{(i)}}] $\delta_0 D$ is analytic and nonvanishing on $\bbD$. 
\item[{\rm{(ii)}}] $\log (\delta_0 D)\in\cap_{p=1}^\infty H^p (\bbD)$ 
\item[{\rm{(iii)}}] For $\f{d\theta}{2\pi}$-a.e.~$e^{i\theta}\in\partial\bbD$ 
with $w(\theta)\neq 0$, 
\begin{equation} \lb{2.18}  
\f{w(\theta)}{w_1(\theta)} = \abs{\delta_0 D (e^{i\theta})}^2  
\end{equation} 
and, in particular, 
\[
\int_{w_1(\theta)\neq 0}  \biggl|\, \log \biggl(\f{w(\theta)}{w_1(\theta)}\biggr) 
\biggr|^p \, \f{d\theta}{2\pi} <\infty 
\]
for all $p\in [1,\infty)$.  
\item[{\rm{(iv)}}] If $\sum_{j=0}^\infty \abs{\alpha_j}^2 <\infty$, then 
\[
(\delta_0 D)(z) = \f{D(z;d\mu)}{D(z;d\mu_1)} 
\]
where $D$ is the Szeg\H{o} function. 
\item[{\rm{(v)}}] If $\varphi_j (z;d\mu_1)$ are the OPUC for $d\mu_1$, then for 
$z\in\bbD$, 
\[
\lim_{n\to\infty} \, \f{\varphi_{n-1}^*(z;d\mu_1)}{\varphi_n^*(z;d\mu)} 
= (\delta_0 D)(z)
\] 
\end{SL} 
\end{theorem}

For a proof, see \cite[Section~2.9]{OPUC1}. The key fact is the calculation 
in $\bbD$ that  
\[
\f{\Real F(z)}{\Real F_1(z)} = \f{\abs{1-\bar\alpha_0 f}^2}{1-\abs{\alpha_0}^2} \, \, 
\f{\abs{1-zf_1}^2}{\abs{1-zf}^2} \,\, \f{1-\abs{z}^2 \abs{f}^2}{1-\abs{f}^2} 
\]
which follows from 
\[
\Real F(z) = \f{1-\abs{z}^2 \abs{f(z)}^2}{1-\abs{f(z)}^2} 
\]
and the Schur algorithm relating $f$ and $f_1$, 
\begin{equation} \lb{2.19}  
zf_1 = \f{f-\alpha_0}{1-\bar\alpha_0 f} 
\end{equation}

\smallskip
One consequence of using $\delta_0 D$ is 

\begin{corollary}\lb{C2.3} Let $d\mu=w(\theta)\f{d\theta}{2\pi}+d\mu_\s$ and 
$d\nu =x(\theta)\f{d\theta}{2\pi} + d\nu_\s$ and suppose that for some $N$ and $k$, 
\[
\alpha_{n+k}(d\mu) =\alpha_n (d\nu) 
\]
for all $n>N$ and that $w(\theta)\neq 0$ for a.e.~$\theta$. Then, $\log (x(\theta)/w(\theta)) 
\in L^1$ and 
\[
\lim_{n\to\infty}\, \f{\|\Phi_n(d\nu)\|^2}{\|\Phi_{n+k}(d\mu)\|^2} = 
\exp\biggl( \int \log \biggl( \f{x(\theta)}{w(\theta)}\biggr) \, \f{d\theta}{2\pi} \biggr) 
\] 
\end{corollary} 

$\delta_0 D$ is also central in the forthcoming paper of Simon-Zlato\v{s} \cite{SZprep}.  

\subsection{Totik's Workshop}\lb{s2.4}

In \cite{Totik}, Totik proved the following: 

\begin{theorem}[Totik \cite{Totik}] \lb{T2.4} Let $d\mu$ be any measure on $\partial\bbD$  
with $\supp (d\mu)=\partial\bbD$. Then there exists a measure $d\nu$ equivalent to $d\mu$  
so that 
\begin{equation} \lb{2.20}  
\lim_{n\to\infty} \, \alpha_n (d\nu)=0 
\end{equation}
\end{theorem} 

This is in a section on Szeg\H{o}'s theorem because Totik's proof uses Szeg\H{o}'s theorem. 
Essentially, the fact that $\sum_{j=0}^\infty \abs{\alpha_j}^2$ doesn't depend on 
$d\mu_\s$ lets one control the a.c.~part of the measure and changes of $\sum_{j=0}^\infty 
\abs{\alpha_j}^2$. By redoing Totik's estimates carefully, one can prove the stronger 
(see \cite[Section~2.10]{OPUC1}): 

\begin{theorem}\lb{T2.5} Let $d\mu$ be any measure on $\partial\bbD$ with $\supp(d\mu) 
=\partial\bbD$. Then there exists a measure $d\nu$ equivalent to $d\mu$ so that for all  
$p>2$, 
\begin{equation} \lb{2.21}  
\sum_{n=0}^\infty \,  \abs{\alpha_n (d\nu)}^p <\infty
\end{equation}
\end{theorem} 

It is easy to extend this to OPRL and there is also a variant for Schr\"odinger 
operators; see Killip-Simon \cite{KSprep}.

\bigskip
\section{The CMV Matrix} \lb{s3} 

One of the most interesting developments in the theory of OPUC in recent years 
is the discovery by Cantero, Moral, and Vel\'azquez \cite{CMV} of a matrix 
realization for multiplication by $z$ on $L^2 (\partial\bbD, d\mu)$ which is 
of finite width (i.e., $\abs{\langle \chi_n, z\chi_m\rangle}=0$ if $\abs{m-n}>k$ 
for some $k$; in this case, $k=2$ to be compared with $k=1$ for OPRL). The obvious 
choice for basis, $\{\varphi_n\}_{n=0}^\infty$, yields a matrix (which \cite{OPUC1} 
calls GGT after Geronimus \cite{Ger44}, Gragg \cite{Gragg82}, and Teplyaev 
\cite{Tep91}) with two defects: If the Szeg\H{o} condition, $\sum_{j=0}^\infty 
\abs{\alpha_j}^2 <\infty$, holds, $\{\varphi_n\}_{n=0}^\infty$ is not a basis 
and $\calG_{k\ell}=\langle \varphi_k, z\varphi_\ell\rangle$ is not unitary. 
In addition, the rows of $\calG$ are infinite, although the columns are finite, 
so $\calG$ is not finite width. 

What CMV discovered is that if $\chi_n$ is obtained by orthonormalizing the 
sequence $1,z,z^{-1}, z^2, z^{-2}, \dots$, we always get a basis 
$\{\chi_n\}_{n=0}^\infty$, in which 
\begin{equation} \lb{3.1}  
\calC_{nm}=\langle \chi_n, z\chi_m\rangle 
\end{equation} 
is five-diagonal. The $\chi$'s can be written in terms of the $\varphi$'s and 
$\varphi^*$ (indeed, $\chi_{2n}=z^{-n}\varphi_{2n}^*$ and $\chi_{2n-1}=
z^{-n+1}\varphi_{2n-1}$) and $\calC$ in terms of the $\alpha$'s. The most 
elegant way of doing this was also found by CMV \cite{CMV}; one can write 
\begin{equation} \lb{3.2}  
\calC=\calL\calM 
\end{equation} 
with 
\begin{equation} \lb{3.3}
\calM = \begin{pmatrix} 
{} & 1 & {} & {} & {} & {} \\
{} & {} & \Theta_1 & {} & {} & {} \\
{} & {} & {} & \Theta_3 & {} & {} \\
{} & {} & {} & {} & \ddots & {} 
\end{pmatrix} 
\qquad 
\calL = \begin{pmatrix} 
{} & \Theta_0 & {} & {} & {} & {} \\
{} & {} & \Theta_2 & {} & {} & {} \\
{} & {} & {} & \Theta_4 & {} & {} \\ 
{} & {} & {} & {} & \ddots & {}
\end{pmatrix} 
\end{equation}  
where the $1$ in $\calM$ is a $1\times 1$ block and all $\Theta$'s are the 
$2\times 2$ block 
\begin{equation} \lb{3.4}
\Theta_j = \begin{pmatrix} 
\bar\alpha_j & \rho_j \\
\rho_j & -\alpha_j \end{pmatrix}
\end{equation}
We let $\calC_0$ denote the CMV matrix for $\alpha_j\equiv 0$. 

The CMV matrix is an analog of the Jacobi matrix for OPRL and it has many 
uses; since \cite{CMV,CMV2} only presented the formalization and a very few 
applications, the section provides numerous new OPUC results based on the CMV 
matrix.  

\subsection{The CMV Matrix and the Szeg\H{o} Function}\lb{s3.1}

If the Szeg\H{o} condition holds, one can define the Szeg\H{o} function 
\begin{equation} \lb{3.5}  
D(z)=\exp \biggl( \int \f{e^{i\theta}+z}{e^{i\theta}-z} \, 
\log (w(\theta))\, \f{d\theta}{4\pi}\biggr)
\end{equation}
One can express $D$ in terms of $\calC$. We use the fact, a special case of 
Lemma~\ref{L3.2} below, that 
\begin{align} 
\sum_{j=0}^\infty \, \abs{\alpha_j}^2 <\infty &\Rightarrow \calC-\calC_0 
\text{ is Hilbert-Schmidt} \lb{3.6}  \\
\sum_{j=0}^\infty\, \abs{\alpha_j}<\infty &\Rightarrow \calC-\calC_0 
\text{ is trace class} \lb{3.7} 
\end{align} 

We also use the fact that if $A$ is trace class, one can define \cite{GK,STI} 
$\det(1+A)$, and if $A$ is Hilbert-Schmidt, $\det_2$ by 
\begin{equation} \lb{3.8}  
\det_2 (1+A)\equiv \det ((1+A)e^{-A}) 
\end{equation}
We also define $w_n$ by 
\begin{equation} \lb{3.9}  
\log (D(z)) =\tfrac12\, w_0 + \sum_{n=1}^\infty z^n w_n 
\end{equation}
so 
\begin{equation} \lb{3.10}  
w_n = \int e^{-in\theta} \log(w(\theta))\, \f{d\theta}{2\pi} 
\end{equation}

Here's the result: 
\begin{theorem}\lb{T3.1} Suppose $\{\alpha_n (d\mu)\}_{n=1}^\infty$ obeys the 
Szeg\H{o} condition  
\begin{equation} \lb{3.11} 
\sum_{n=0}^\infty \, \abs{\alpha_n}^2 <\infty 
\end{equation} 
Then the Szeg\H{o} function, $D$, obeys for $z\in\bbD$,
\begin{equation} \lb{3.12} 
D(0)D(z)^{-1} = \det_2\biggl( \f{(1-z\bar\calC)}{(1-z\bar\calC_0)}\biggr) e^{+zw_1} 
\end{equation}
where 
\begin{equation} \lb{3.13} 
w_1 = \alpha_0 - \sum_{n=1}^\infty \alpha_n \bar\alpha_{n-1} 
\end{equation}
If 
\begin{equation} \lb{3.14} 
\sum_{n=0}^\infty \abs{\alpha_n}<\infty
\end{equation}
then 
\begin{equation} \lb{3.15} 
D(0)D(z)^{-1} = \det\biggl( \f{(1-z\bar\calC)}{(1-z\bar\calC_0)}\biggr) 
\end{equation}

The coefficients $w_n$ of \eqref{3.9} are given by 
\begin{equation} \lb{3.16} 
w_n = \f{\, \ol{\tr(\calC^n - \calC_0^n)}\,}{n}
\end{equation}
for all $n\geq 1$ if \eqref{3.14} holds and for $n\geq 2$ if \eqref{3.11} holds. In all cases, 
one has 
\begin{equation} \lb{3.17} 
w_n = \sum_{j=0}^\infty \f{\ol{(\calC^n)}_{jj}}{n}
\end{equation}
\end{theorem}

{\it Remark.} $\bar\calC$ is the matrix $(\bar\calC)_{k\ell} = \ol{(\calC_{k\ell})}$. 

\smallskip
The proof (given in \cite[Section~4.2]{OPUC1}) is simple: by \eqref{4.10} below, $\Phi_n$ can 
be written as a determinant of a cutoff CMV matrix, which gives a formula for $\varphi_n^*$. 
Since $\varphi_n^*\to D^{-1}$, the cutoff matrices converge in Hilbert-Schmidt and 
trace norm and since $\det$/$\det_2$ are continuous, one can take limits of the 
finite formulae.

\subsection{CMV Matrices and Spectral Analysis} \lb{s3.2}

The results in this subsection are joint with Leonid Golinskii. The CMV matrix provides 
a powerful tool for the comparison of properties of two measures $d\mu$, $d\nu$ on 
$\partial\bbD$ if we know something about $\alpha_n (d\nu)$ as a perturbation of 
$\alpha_n (d\mu)$. Of course, this idea is standard in OPRL and Schr\"odinger 
operators. For example, Krein \cite{AK} proved a theorem of Stieltjes \cite{Stie} that 
$\supp (d\mu)$ has a single non-isolated point $\lambda$ if and only if the Jacobi 
parameters $a_n\to 0$ and $b_n\to\lambda$ by noting both statements are equivalent to 
$J-\lambda 1$ being constant. Prior to results in this section, many results were proven 
using the GGT representation, but typically, they required $\liminf_{n\to\infty} 
\abs{\alpha_n}>0$ to handle the infinite rows. 

Throughout this section, we let $d\mu$ (resp.~$d\nu$) have Verblunsky coefficient 
$\alpha_n$ (resp.~$\beta_n$) and we define $\rho_n=(1-\abs{\alpha_n}^2)^{1/2}$, 
$\sigma_n =(1-\abs{\beta_n})^{1/2}$. An easy estimate using the $\calL\calM$ 
factorization  shows with $\|\cdot\|_p$ the $\calI_p$ trace ideal norm 
\cite{GK,STI}:

\begin{lemma}\lb{L3.2} There exists a universal constant $C$ so that for all 
$1\leq p\leq\infty$, 
\begin{equation} \lb{3.18}  
\|\calC(d\mu)-\calC(d\nu)\|_p\leq C \biggl(\, \sum_{n=0}^\infty \, 
\abs{\alpha_n-\beta_n}^p + \abs{\rho_n -\sigma_n}^p\biggr)^{1/p} 
\end{equation}
\end{lemma} 

{\it Remark.} One can take $C=6$. For $p=\infty$, the right side of \eqref{3.18} 
is interpreted as $\sup_n (\max (\abs{\alpha_n-\beta_n}, \abs{\rho_n-\sigma_n}))$. 

\smallskip
This result allows one to translate the ideas of Simon-Spencer \cite{S208} to a 
new proof of the following result of Rakhmanov \cite{Rakh83} (sometimes called 
Rakhmanov's lemma): 

\begin{theorem}\lb{T3.3} If $\limsup \abs{\alpha_n}=1$, $d\mu$ is purely singular. 
\end{theorem} 

\begin{proof}[Sketch] Pick a subsequence $n_j$ so 
\begin{equation} \lb{3.19}  
\sum_{j=0}^\infty \, (1-\abs{\alpha_{n_j}})^{1/2} <\infty 
\end{equation} 
Let $\beta_k =\alpha_k$ if $k\neq n_j$ and $\beta_k=\alpha_k/\abs{\alpha_k}$ if 
$k=n_j$. There is a limiting unitary $\ti\calC$ with those values of $\beta$. It  
is a direct sum of finite rank matrices since $\abs{\beta_{n_j}}=1$ forces $\calL$ 
or $\calM$ to have some zero matrix elements. Thus $\ti\calC$ has no a.c.~spectrum. 

By \eqref{3.19} and \eqref{3.18}, $\calC-\ti\calC$ is trace class, so by the the 
Kato-Birman theorem for unitaries \cite{BK62}, $\calC$ has simply a.c.~spectrum. 
\end{proof} 

Golinskii-Nevai \cite{GN} already remarked that Rakhmanov's lemma is an analog of 
\cite{S208}. For the next pair of results, the special case $\lambda_n\equiv 1$ 
are analogs of extended results of Weyl and Kato-Birman but for OPUC are new 
even in this case with the generality we have. 

\begin{theorem}\lb{T3.4} Suppose $\{\lambda_n\}_{n=0}^\infty \in\partial\bbD^\infty$, 
$\{\alpha_n\}_{n=0}^\infty,\{b_n\}_{n=0}^\infty \in\bbD^\infty$ and 
\begin{alignat*}{2} 
&\text{\rm{(i)}} \qquad &&\beta_n \lambda_n -\alpha_n \to 0 \\
&\text{\rm{(ii)}} \qquad && \lambda_{n-1} \bar\lambda_n\to 1 
\end{alignat*}
Then the derived sets of $\supp(d\mu)$ and $\supp(d\nu)$ are equal, that is, up to 
a discrete set, $\supp(d\mu)$ and $\supp(d\nu)$ are equal. 
\end{theorem} 

\begin{theorem}\lb{T3.5} Suppose $\{\lambda_n\}_{n=0}^\infty\in\partial\bbD^\infty$ 
and $\alpha_n,\beta_n$ are the Verblunsky coefficients of $d\mu =w(\theta) 
\f{d\theta}{2\pi} +d\mu_\s$ and $d\nu = f(\theta) \f{d\theta}{2\pi} + d\nu_\s$. 
Suppose that 
\[
\sum_{j=0}^\infty \, \abs{\lambda_j \alpha_j-\beta_j} + 
\abs{\lambda_{j+1} \bar\lambda_j-1} <\infty 
\]
Then $\{\theta\mid w(\theta)\neq 0\}=\{\theta\mid f(\theta)\neq 0\}$ {\rm{(}}up to 
sets of $d\theta/2\pi$ measure $0${\rm{)}}. 
\end{theorem} 
 
The proofs (see \cite[Section~4.3]{OPUC1}) combine the estimates of Lemma~\ref{L3.2} 
and the fact that conjugation of CMV matrices with diagonal matrices can be realized as 
phase changes. That $\supp (d\mu)=\partial\bbD$ if $\abs{\alpha_j}\to 0$ (special case 
of Theorem~\ref{T3.4}) is due to Geronimus \cite{GBk1}. Other special cases can 
be found in \cite{BRLL,Gol2000A}. 

\cite[Section~4.3]{OPUC1} also has results that use trial functions and CMV matrices. 
Trial functions are easier to use for unitary operators than for selfadjoint ones since 
linear variational principles for selfadjoint operators only work at the ends of the 
spectrum. But because $\partial\bbD$ is curved, linear variational principles work 
at any point in $\partial\bbD$. For example, $(\theta_0-\veps, \theta_0 +\veps)\cap 
\supp(d\mu)=\emptyset$ if and only if  
\[
\Real (e^{-i\theta_0} \langle \psi, (e^{i\theta_0}- \calC)\psi\rangle) \geq 
2\sin^2 \biggl( \f{\veps}{2}\biggr) \|\psi\|^2 
\] 
for all $\psi$. Typical of the results one can prove using trial functions is: 

\begin{theorem}\lb{T3.6} Suppose there exists $N_j\to\infty$ and $k_j$ so 
\[
\f{1}{N_j} \sum_{\ell=1}^{N_j}\, \abs{\alpha_{k_j+\ell}}^2\to 0
\] 
Then $\supp (d\mu)=\partial\bbD$. 
\end{theorem}

\subsection{CMV Matrices and the Density of Zeros} \lb{s3.3}

A fundamental object of previous study is the density of zeros, $d\nu_n (z;d\mu)$, 
defined to give weight $k/n$ to a zero of $\Phi_n (z;d\mu)$ of multiplicity $k$. 
One is interested in its limit or limit points as $n\to\infty$. A basic difference 
from OPRL is that for OPRL, any limit point is supported on $\supp(d\mu)$, while 
limits of $d\nu_n$ need not be supported on $\partial\bbD$. Indeed, for $d\mu = 
d\theta/2\pi$, $d\nu_n$ is a delta mass at $z=0$ and \cite{SimTot} have found $d\mu$'s 
for which the limit points of $d\nu_n$ are all measures on $\bar\bbD$!  

As suggested by consideration of the ``density of states" for Schr\"odinger operators 
and OPRL (see \cite{P73,S149}), moments of the density of zeros are related to traces 
of powers of a truncated CMV matrix. Define $\calC^{(n)}$ to be the matrix obtained 
from the topmost $n$ rows and leftmost $n$ columns of $\calC$. Moreover, let $d\gamma_n$ 
be the Ces\`aro mean of $\abs{\varphi_j}^2\, d\mu$, that is, 
\begin{equation} \lb{3.20}  
d\gamma_n (\theta) =\f{1}{n} \sum_{j=0}^{n-1} \, \abs{\varphi_j (e^{i\theta},d\mu)}^2\,   
d\mu (\theta)
\end{equation} 
Then: 

\begin{theorem}\lb{T3.7} For any $k\geq 0$, 
\begin{equation} \lb{3.21}  
\int z^k \, d\nu_n (z)=\f{1}{n} \, \tr ((\calC^{(n)})^k) 
\end{equation}
Moreover, 
\begin{equation} \lb{3.22}  
\lim_{n\to\infty}\, \biggl[\biggl( \int z^k \, d\nu_n(z)\biggr) - \biggl( \int z^k\, 
d\gamma_n (z)\biggr)\biggr] =0
\end{equation} 
\end{theorem} 

\begin{proof}[Sketch] (For details, see \cite[Section~8.2]{OPUC1}.) We'll see in 
Theorem~\ref{T4.5} that the eigenvalues of $\calC^{(n)}$ (counting geometric multiplicity) 
are the zeros of $\Phi_n (z;d\mu)$ from which \eqref{3.21} is immediate. 

Under the CMV representation, $\delta_j$ corresponds to $z^\ell \varphi_j$ or 
$z^\ell \varphi_j^*$ for suitable $\ell$ (see the discussion after \eqref{3.1}) so 
\[
(\calC^k)_{jj}=\int e^{ik\theta} \abs{\varphi_j (e^{i\theta})}^2\, d\mu(\theta) 
\]
and thus  
\begin{equation} \lb{3.23}  
\int z^k\, d\gamma_n(z) = \f{1}{n}\sum_{j=0}^{n-1} (\calC^k)_{jj} 
\end{equation}

If $\ell <n-2k$, 
\[
([\calC^{(n)}]^k)_{\ell\ell} = (\calC^k)_{\ell\ell} 
\]
so that \eqref{3.22} follows from \eqref{3.21} and \eqref{3.23}. 
\end{proof} 

From \eqref{3.21} and \eqref{3.18}, we immediately get 

\begin{corollary}\lb{C3.7A} If $\lim_{N\to\infty} \f{1}{N}\sum_{j=0}^{N-1} 
\abs{\alpha_j-\beta_j}\to 0$, then for any $k$, 
\begin{equation} \lb{3.23A} 
\lim_{N\to\infty} \, \int z^k [d\nu_N(z;\{\alpha_j\}_{j=0}^\infty) - 
d\nu_N (z;\{\beta_j\}_{j=0}^\infty)]=0
\end{equation} 
\end{corollary} 

One application of this is to a partially alternative proof of a theorem of 
Mhaskar-Saff \cite{MhS1}. They start with an easy argument that uses a theorem of 
Nevai-Totik \cite{NT89} and the fact that $(-1)^{n+1} \bar\alpha_{n-1}$ is the product 
of zeros of $\Phi_n(z)$ to prove 

\begin{lemma}\lb{L3.8} Let 
\begin{equation} \lb{3.24} 
A=\limsup \, \abs{\alpha_n}^{1/n} 
\end{equation}
and pick $n_j$ so 
\begin{equation} \lb{3.25} 
\abs{\alpha_{n_j-1}}^{1/n_j-1}\to A
\end{equation}
Then any limit points of $d\nu_{n_j}$ is supported on $\{z\mid \abs{z}=A\}$. 
\end{lemma} 

They then use potential theory to prove the following, which can be proven instead 
using the CMV matrix: 

\begin{theorem}[Mhaskar-Saff \cite{MhS1}] \lb{T3.9} Suppose \eqref{3.24} and \eqref{3.25} 
hold and that either $A<1$ or 
\begin{equation} \lb{3.26} 
\lim_{n\to\infty}\, \f{1}{N} \sum_{j=0}^{n-1}\, \abs{\alpha_j}=0
\end{equation}
Then $d\nu_{n_j}$ converges weakly to the uniform measure on the set $\{z\mid\abs{z}=A\}$. 
\end{theorem} 

\begin{proof}[Sketch of New Proof] Since $d\theta/2\pi$ is the unique measure with 
$\int z^k \f{d\theta}{2\pi}=\delta_{k0}$ for $k\geq 0$, it suffices to show that 
for $k\geq 1$, 
\[
\int z^k\, d\nu_{n_j} \to 0 
\]
This is immediate from Corollary~\ref{C3.7A} and the fact that $\int z^k d\ti\nu_n=0$ 
if $\ti\nu_n$ is the zero's measure for $d\theta/2\pi$. 
\end{proof}

\subsection{CMV and Wave operators}\lb{s3.4}

In \cite[Section~10.7]{OPUC2}, we prove the following: 

\begin{theorem}\lb{T3.10} Suppose $\sum_{n=0}^\infty \abs{\alpha_n}^2<\infty$. Let 
$\calC$ be the CMV matrix for $\{\alpha_n\}_{n=0}^\infty$ and $\calC_0$ the CMV matrix 
for $\alpha_j\equiv 0$. Then 
\[
\slim_{n\to\pm\infty}\, \calC^n \calC_0^{-n} =\Omega^\pm 
\]
exists and its range is $\chi_S(\calC)$ where $S$ is a set with $d\mu_\s(S)=0$, 
$\abs{\partial\bbD\backslash S}=0$. 
\end{theorem} 

The proof depends on finding an explicit formula for $\Omega^\pm$ (in terms of $D(z)$, 
the Szeg\H{o} function); equivalently, from the fact that in a suitable sense, $\calC$ 
has no dispersion. The surprise is that one only needs $\sum_{n=0}^\infty \abs{\alpha_n}^2 
<\infty$, not $\sum_{n=0}^\infty \abs{\alpha_n}<\infty$. Some insight can be obtained 
from the formulae Geronimus \cite{Ger46} found mapping to a Jacobi matrix when the 
$\alpha$'s are real. The corresponding $a$'s and $b$'s have the form $c_{n+1}-c_n +d_n$ 
where $d_n\in\ell^1$ and $c_n\in\ell^2$, so there are expected to be modified wave 
operators with finite modifications since $c_{n+1}-c_n$ is conditionally summable. 

Simultaneous with our discovery of Theorem~\ref{T3.10}, Denisov \cite{DenGAFA} found a 
similar result for Dirac operators.

\subsection{The Resolvent of the CMV Matrix}\lb{s3.5}

I have found an explicit formula for the resolvent of the CMV matrix $(\calC-z)_{k\ell}^{-1}$ 
when $z\in\bbD$ (and for some suitable limits as $z\to\partial\bbD$), not unrelated to a 
formula for the resolvent of the GGT matrix found by Geronimo-Teplyaev \cite{GTep} 
(see also \cite{GJo2,GJo1}).  

Just as the CMV basis, $\chi_n$, is the result of applying Gram-Schmidt to orthonormalize 
$\{1,z,z^{-1}, z^2, z^{-2}, \dots\}$, the alternate CMV basis, $x_n$, is what we get by 
orthonormalizing $\{1,z^{-1}, z,z^{-2}, z^2, \dots\}$. (One can show $\ti\calC=\langle 
x_{\boldsymbol{\cdot}}, zx_{\boldsymbol{\cdot}}\rangle =\calM\calL$.) Similarly, let 
$y_n, \Upsilon_n$ be the CMV and alternate CMV bases associated to $(\psi_n, -\psi_n^*)$. 
Define 
\begin{align} 
p_n &= y_n + F(z) x_n \lb{3.27} \\
\pi_n & = \Upsilon_n + F(z) \chi_n \lb{3.28}
\end{align}
Then 

\begin{theorem}\lb{T3.11}  We have that for $z\in\bbD$, 
\begin{equation} \lb{3.29}
[(\calC-z)^{-1}]_{k\ell} = 
\begin{cases} (2z)^{-1}\chi_\ell(z) p_k (z) & k>\ell \text{ or } k=\ell =2n-1 \\ 
(2z)^{-1} \pi_\ell (z) x_k (z) & \ell >k \text{ or } k=\ell=2n \end{cases}
\end{equation}
\end{theorem}  

This is proven in \cite[Section~4.4]{OPUC1}. It can be used to prove Khrushchev's 
formula \cite{Kh2000} that the Schur function for $\abs{\varphi_n}^2\,d\mu$ is 
$\varphi_n (\varphi_n^*)^{-1} f(z;\{\alpha_{n+j}\}_{j=0}^\infty)$; see 
\cite[Section~9.2]{OPUC2}.

\subsection{Rank Two Perturbations and CMV Matrices}\lb{s3.6}

We have uncovered some remarkably simple formulae for finite rank perturbations of 
unitaries. If $U$ and $V$ are unitary so $U\varphi=V\varphi$ for $\varphi\in 
\ran (1-P)$ where $P$ is a finite-dimensional orthogonal projection, then there is 
a unitary $\Lambda =P\calH\to P\calH$ so that 
\begin{equation} \lb{3.29A} 
V=U(1-P)+ U\Lambda P
\end{equation} 
For $z\in\bbD$, define $G_0(z)$, $G(z)$, $g_0(z)$, $g(z)$ mapping $P\calH$ to $P\calH$ by 
\begin{align} 
G(z) &= P \biggl[ \f{V+z}{V-z}\biggr] P \lb{3.30} \\
G_0(z) &= P\biggl[ \f{U+z}{U-z}\biggr] P \lb{3.31} 
\end{align}
\begin{equation} \lb{3.32}
G(z)=\f{1+zg(z)}{1-zg(z)} \qquad G_0(z) = \f{1+zg_0(z)}{1-zg_0(z)}
\end{equation} 
As operators on $P\calH$, $\|g(z)\|<1$, $\|g_0(z)\|<1$ on $\bbD$. A direct calculation 
(see \cite[Section~4.5]{OPUC1}) proves that 
\begin{equation} \lb{3.33} 
g(z)=\Lambda^{-1} g_0(z)
\end{equation}
This can be used to provide, via a rank two decoupling of a CMV matrix (change a 
$\Theta(\alpha)$ to $\left(\begin{smallmatrix} -1 & 0 \\ 0 & 1\end{smallmatrix}\right)$),  
new proofs of Geronimus' theorem and of Khrushchev's formula; see \cite[Section~4.5]{OPUC1}.

\subsection{Extended and Periodized CMV Matrices}\lb{s3.7}

The CMV matrix is defined on $\ell^2 (\{0,1,\dots\})$. It is natural to define an extended 
CMV matrix associated to $\{\alpha_j\}_{j=-\infty}^\infty$ on $\ell^2 (\bbZ)$ by extending 
$\calL$ and $\calM$ to $\ti\calL$ and $\ti\calM$ on $\ell^2 (\bbZ)$ as direct sums of 
$\Theta$'s and letting $\calE=\ti\calL\ti\calM$. 

This is an analog of whole-line discrete Schr\"odinger operators. It is useful in the study 
of OPUC with ergodic Verblunsky coefficients as well as a natural object in its own 
right. \cite{OPUC1,OPUC2} have numerous results about this subject introduced here for the 
first time. 

If $\{\alpha_j\}_{j=-\infty}^\infty$ is periodic of period $p$, $\calE$ commutes with 
translations and so is a direct integral of $p\times p$ periodized CMV matrices depending 
on $\beta\in\partial\bbD$: essentially to restrictions of $\calE$ to sequences in $\ell^\infty$ 
with $u_{n+kp}=\beta^k u_n$. In \cite[Section~12.1]{OPUC2}, these are linked to Floquet 
theory and to the discriminant, as discussed below in Section~\ref{s5.1}.

\bigskip
\section{Miscellaneous Results} \lb{s4} 

In this section, we discuss a number of results that don't fit into the themes of 
the prior sections and don't involve explicit classes of Verblunsky coefficients, 
the subject of the final two sections.

\subsection{Jitomirskaya-Last Inequalities}\lb{s4.1} 

In a fundamental paper intended to understand the subordinacy results of Gilbert-Pearson 
\cite{GP} and extend the theory to understand Hausdorff dimensionality, Jitomirskaya-Last 
\cite{JL0,JL1} proved some basic inequalities about singularities of the $m$-function 
as energy approaches the spectrum. 

In \cite[Section~10.8]{OPUC2}, we prove an analog of their result for OPUC. First, we 
need some notation. $\psi_{\boldsymbol{\cdot}}$ denotes the second polynomial, that 
is, the OPUC with sign flipped $\alpha_j$'s. For $x\in [0,\infty)$, let $[x]$ be 
the integral part of $x$ and define for a sequence $a$: 
\begin{equation} \lb{4.1}
\|a\|_x^2 =\sum_{j=0}^{[x]}\, \abs{a_j}^2 + (x-[x]) \abs{a_{j+1}}^2 
\end{equation}
We prove 

\begin{theorem}\lb{T4.1} For $z\in\partial\bbD$ and $r\in [0,1)$, define $x(r)$ to be 
the unique solution of 
\begin{equation}\lb{4.2} 
(1-r)\|\varphi_{\boldsymbol{\cdot}}(z)\|_{x(r)} \|\psi_{\boldsymbol{\cdot}}(z)\|_{x(r)} 
=\sqrt{2}
\end{equation}
Then 
\begin{equation}\lb{4.3} 
A^{-1} \biggl[ \f{\|\psi_{\boldsymbol{\cdot}}(z)\|_{x(r)}}{\|\varphi_{\boldsymbol{\cdot}}(z)\|_{x(r)}} 
\biggr] \leq \abs{F(rz)} \leq A \biggl[ \f{\|\psi_{\boldsymbol{\cdot}}(z)\|_{x(r)}} 
{\|\varphi_{\boldsymbol{\cdot}}(z)\|_{x(r)}}\biggr]
\end{equation}
where $A$ is a universal constant in $(1,\infty)$. 
\end{theorem}  

{\it Remark.} One can take $A=6.65$; no attempt was made to optimize $A$. 

\smallskip
This result allows one to extend the Gilbert-Pearson subordinacy theory \cite{GP} 
to OPUC. Such an extension was accomplished by Golinskii-Nevai \cite{GN} under an 
extra assumption that 
\begin{equation} \lb{4.4}
\limsup \, \abs{\alpha_n} <1
\end{equation}
We do not need this assumption, but the reason is subtle as we now explain.  

Solutions of \eqref{1.4} and its ${}^*{}$ viewed as an equation for 
$\binom{\varphi}{\varphi^*}$ are given by a transfer matrix 
\begin{equation} \lb{4.4a}
T_n(z) = A(\alpha_{n-1},z) A(\alpha_{n-2}, z) \dots A(\alpha_0,z)
\end{equation}
where $\rho = (1-\abs{\alpha}^2)^{1/2}$ and 
\begin{equation} \lb{4.4b}
A(\alpha,z)=\rho^{-1} \begin{pmatrix} z & -\bar\alpha \\
-\alpha z & 1 \end{pmatrix}
\end{equation}

In the discrete Schr\"odinger case, the transfer matrix is a product of $A(v,e)= 
\left(\begin{smallmatrix} e-v & -1 \\ 1 & 0\end{smallmatrix}\right)$. A key role 
in the proof in \cite{JL1} is that $A(v,e')-A(v,e)$ depends only on $e$ and $e'$ 
and not on $v$. For OPUC, the $A$ has the form \eqref{4.4b}. \cite{GN} requires 
\eqref{4.4} because $A(\alpha,z)-A(\alpha,z')$ has a $\rho^{-1}$ divergence, and 
\eqref{4.4} controlled that. The key to avoiding \eqref{4.4} is to note that 
\[
A(\alpha,z)-A(\alpha,w) =(1-z^{-1} w) A(\alpha,z)P
\]
where $P=\left(\begin{smallmatrix}1 & 0 \\ 0 & 0\end{smallmatrix}\right)$.

\subsection{Isolated Pure Points}\lb{s4.2}

Part of this section is joint work with S.~Denisov. These results extend beyond 
the unit circle. We'll be interested in general measures on $\bbC$ with nontrivial 
probability measures 
\begin{equation} \lb{4.5}
\int \abs{z}^j \, d\mu(z) <\infty
\end{equation}
for all $j=0,1,2,\dots$. In that case, one can define monic orthogonal polynomials 
$\Phi_n(z)$, $n=0,1,2,\dots$. Recall the following theorem of Fej\'er \cite{Fej}: 

\begin{theorem}[Fej\'er \cite{Fej}]\lb{T4.2} All the zeros of $\Phi_n$ lie in the  
convex hull of $\supp(d\mu)$. 
\end{theorem} 

We remark that this theorem has an operator theoretic interpretation. If $M_z$ is 
the operator of multiplication by $z$ on $L^2 (\bbC,d\mu)$, and if $P_n$ is the projection 
onto the span of $\{z^j\}_{j=0}^{n-1}$, then we'll see \eqref{4.8} that the eigenvalues 
of $P_nM_zP_n$ are precisely the zeros of $\Phi_n$. If $\eta (\dott)$ denotes 
numerical range, $\eta (M_z)$ is the convex hull of $\supp(d\mu)$, so Fej\'er's 
theorem follows from $\eta (P_nM_zP_n)\subseteq \eta(M_z)$ and the fact that 
eigenvalues lie in the numerical range. 

\cite[Section~1.7]{OPUC1} contains the following result I proved with Denisov: 

\begin{theorem}\lb{T4.3} Let $\mu$ obey \eqref{4.5} and suppose $z_0$ is an isolated 
point of $\supp(d\mu)$. Define $\Gamma =\supp(d\mu)\backslash\{z_0\}$ and $\ch (\Gamma)$, 
the convex hull of $\Gamma$. Suppose $\delta\equiv\dist (z_0,\ch(\Gamma))>0$. Then 
$\Phi_n$ has at most one zero in $\{z\mid\abs{z-z_0}< \delta/3\}$. 
\end{theorem} 

{\it Remarks.} 1. In case $\supp(d\mu)\subset\partial\bbD$, any isolated point has 
$\delta >0$. Indeed, if $d=\dist (z_0,\Gamma)$, $\delta  \geq d^2/2$ and so, 
Theorem~\ref{T4.3} says that there is at most one zero in the circle of radius 
$d^2/6$. 

\smallskip
2. If $d\mu$ is a measure on $[-1,-\f12]\cup\{0\}\cup [\f12,1]$ and symmetric under 
$x$, and $\mu(\{0\})>0$, it can be easily shown that $P_{2n}(x)$ has two zeros near 
$0$ for $n$ large. Thus, for a result like Theorem~\ref{T4.3}, it is not enough that 
$z_0$ be an isolated point of $\supp(d\mu)$; note in this example that $0$ is in the 
convex hull of $\supp(d\mu)\backslash \{0\}$.  

\smallskip
The other side of this picture is the following result proven in \cite[Section~8.1]{OPUC1} 
using potential theoretic ideas of the sort exposed in \cite{SaffTot,SST}: 

\begin{theorem} \lb{T4.4} Let $\mu$ be a nontrivial probability measure on $\partial 
\bbD$ and let $z_0$ be an isolated point of $\supp(d\mu)$. Then there exist $C>0$, 
$a>0$, and a zero $z_n$ of $\Phi_n(z;d\mu)$ so that 
\begin{equation} \lb{4.6}
\abs{z_n-z_0} \leq Ce^{-a\abs{n}}
\end{equation}
\end{theorem} 

There is an explicit formula for $a$ in terms of the equilibrium potential for 
$\supp(d\mu)$ at $z_0$. The pair of theorems in this section shows that any isolated 
mass point, $z_0$, of $d\mu$ on $\partial\bbD$ has exactly one zero near $z_0$ for 
$n$ large.

\subsection{Determinant Theorem} \lb{s4.3}

It is a well-known fact that if $J^{(n)}$ is the $n\times n$ truncated Jacobi matrix 
and $P_n$ the monic polynomial associated to $J$, then 
\begin{equation} \lb{4.7}
P_n(x)=\det (x-J^{(n)})
\end{equation}
The usual proofs of \eqref{4.7} use the selfadjointness of $J^{(n)}$ but there is a 
generalization to OPs for measures on $\bbC$: 

\begin{theorem}\lb{T4.5} Let $d\mu$ be a measure on $\bbC$ obeying \eqref{4.5}. 
Let $P_n$ be the projection onto the span of $\{z^j\}_{j=0}^{n-1}$, $M_z$ be 
multiplication by $z$, and $M^{(n)}=P_nM_zP_n$. Then 
\begin{equation} \lb{4.8}
\Phi_n(z) =\det (z-M^{(n)})
\end{equation}
\end{theorem} 

\begin{proof}[Sketch] Suppose $z_0$ is an eigenvalue of $M^{(n)}$. Then there exists 
$Q$, a polynomial of degree at most $n-1$, so $P_n(z-z_0) Q(z)=0$. Since $\Phi_n$ is 
up to a constant, the only polynomial, $S$, of degree $n$ with $P_n(S)=0$, we see 
\begin{equation} \lb{4.9}
(z-z_0)Q(z)=c\Phi_n(z)
\end{equation}
It follows that $\Phi_n(z_0)=0$, and conversely, if $\Phi_n(z_0)=0$, $\Phi_n(z)/(z-z_0) 
\equiv Q$ provides an eigenfunction. Thus, the eigenvalues of $M^{(n)}$ are exactly 
the zeros of $\Phi_n$. This proves \eqref{4.8} if $\Phi_n$ has simple zeros. 
In general, by perturbing $d\mu$, we can get $\Phi_n$ as a limit of other $\Phi_n$'s 
with simple zeros. 
\end{proof} 

In the case of $\partial\bbD$, $z^\ell$ is unitary on $L^2 (\partial\bbD,d\mu)$, so 
$P_n$ in defining $M^{(n)}$ can be replaced by the projection onto the span of 
$\{z^{j+\ell}\}_{j=0}^{n-1}$ for any $\ell$, in particular, the span onto the first 
$n$ of $1,z,z^{-1},z^2, \dots$, so 

\begin{corollary}\lb{C4.6} If $\calC^{(n)}$ is the truncated $n\times n$ CMV matrix, then 
\begin{equation} \lb{4.10}
\Phi_n(z)=\det (z-\calC^{(n)})
\end{equation}
\end{corollary}

\subsection{Geronimus' Theorem and Taylor Series} \lb{s4.4}

Given a Schur function, that is, $f$ mapping $\bbD$ to $\bbD$ analytically, one 
defines $\gamma_0$ and $f_1$ by 
\begin{equation} \lb{4.11}
f(z)= \f{\gamma_0 + zf_1(z)}{1+\bar\gamma_0 zf_1(z)} 
\end{equation} 
so $\gamma_0 =f(0)$ and $f_1$ is either a new Schur function or a constant in $\partial\bbD$. 
The later combines the fact that $\omega\to (\gamma_0+\omega)/(1+\bar\gamma_0 \omega)$ is 
a bijection of $\bbD$ to $\bbD$ and the Schur lemma that if $g$ is a Schur function with 
$g(0)=0$, then $g(z)z^{-1}$ is also a Schur function. If one iterates, one gets either a 
finite sequence $\gamma_0, \dots, \gamma_{n-1}\in\bbD^n$ and $\gamma_n \in\partial\bbD$ 
or an infinite sequence $\{\gamma_j\}_{j=0}^\infty \in\bbD^\infty$. It is a theorem 
of Schur that this sets up a one-one correspondence between the Schur functions and such 
$\gamma$-sequences. The finite sequences correspond to finite Blaschke products. 

In 1944, Geronimus proved 

\begin{theorem}[Geronimus' Theorem \cite{Ger44}]\lb{T4.7} Let $d\mu$ be a nontrivial 
probability measure on $\partial\bbD$ with Verblunsky coefficients $\{\alpha_j\}_{j=0}^\infty$. 
Let $f$ be the Schur function associated to $d\mu$ by \eqref{2.14}/\eqref{2.16} and let  
$\{\gamma_n\}_{n=0}^\infty$ be its Schur parameters. Then 
\begin{equation} \lb{4.11A}
\gamma_n =\alpha_n
\end{equation} 
\end{theorem} 

\cite{OPUC1} has several new proofs of this theorem (see \cite{Gol93,PN,K91} for other 
proofs, some of them also discussed in \cite{OPUC1}). We want to describe here one proof 
that is really elementary and should have been found in 1935! Indeed, it is obvious to 
anyone who knows Schur's paper \cite{Sc1917} and Verblunsky \cite{V35} --- but apparently 
Verblunsky didn't absorb that part of Schur's work, and Verblunsky's paper seems to have 
been widely unknown and unappreciated! 

This new proof depends on writing the Taylor coefficients of $F(z)$ in terms of the 
$\alpha$'s and the $\gamma$'s. Since 
\[
\f{e^{i\theta}+z}{e^{i\theta}-z} = 1+2\sum_{n=1}^\infty e^{-in\theta} z^n 
\]
we have 
\begin{equation} \lb{4.12}
F(z)=1+2 \sum_{n=1}^\infty c_n z^n 
\end{equation}
with $c_n$ given by 
\begin{equation} \lb{4.13}
c_n =\int e^{-in\theta} \, d\mu(\theta) 
\end{equation}

Define $s_n(f)$ by $f(z)=\sum_{n=0}^\infty s_n(f) z^n$. Then Schur \cite{Sc1917} noted 
that $(1+\bar\gamma_0 zf_1)f=\gamma_0 + zf_1$ implies 
\[
s_n(f) = (1-\abs{\gamma_0}^2) s_{n-1} (f_1) -\bar\gamma_0 \sum_{j=1}^n s_j(f) 
s_{n-1-j} (f_1) 
\]
so that, by induction,  
\[
s_n(f) =\prod_{j=0}^{n-1} (1-\abs{\gamma_j}^2) \gamma_n + r_n (\gamma_0, \bar\gamma_0,  
\dots, \gamma_{n-1}, \bar\gamma_{n-1}) 
\]
with $r_n$ a polynomial. This formula is in Schur \cite{Sc1917}. Since $F(z)=1+2 
\sum_{n=1}^\infty (zf)^n$, we find that $c_n=s_{n-1}(f)+ \text{ polynomial}$ in $(s_0(f), 
\dots, s_{n-1}(f))$, and thus 
\begin{equation} \lb{4.14}
c_n(f) =\prod_{j=0}^{n-2} (1-\abs{\gamma_j}^2) \gamma_{n-1} + \ti r_{n-1} 
(\gamma_0, \bar\gamma_0, \dots, \gamma_{n-2}, \bar\gamma_{n-2})
\end{equation} 
for a suitable polynomial $\ti r_{n-1}$. 

On the other hand, Verblunsky \cite{V35} had the formula relating his parameters and 
$c_n (f)$: 
\begin{equation} \lb{4.15}
c_n(f) =\prod_{j=0}^{n-2} (1-\abs{\alpha_n}^2) \alpha_{n-1} + \ti q_{n-1} 
(\alpha_0, \bar\alpha_0, \dots, \alpha_{n-2}, \bar\alpha_{n-2}) 
\end{equation} 
For Verblunsky, \eqref{4.15} was actually the definition of $\alpha_{n-1}$, that is, 
he showed (as did Akhiezer-Krein \cite{AK36}) that, given $c_0, \dots, c_{n-1}$, 
the set of allowed $c_n$'s for a positive Toeplitz determinant is a circle of 
radius inductively given by $\prod_{j=0}^{n-1} (1-\abs{\alpha_j}^2)$, which led him 
to define parameters $\alpha_{n-1}$. 

On the other hand, it is a few lines to go from the Szeg\H{o} recursion \eqref{1.4} to 
\eqref{4.15}. For we note that 
\[
\int \Phi_{n+1}(z) \, d\mu(z) = \langle 1, \Phi_{n+1}\rangle =0 
\]
while 
\[
\langle 1, \Phi_n^*\rangle = \langle \Phi_n,z^n\rangle = 
\langle \Phi_n, \Phi_n\rangle = \prod_{j=1}^{n-1} (1-\abs{\alpha_j}^2) 
\]
by \eqref{1.7}. Thus 
\begin{equation} \lb{4.16}
\langle 1, z\Phi_n\rangle =\bar\alpha_n \prod_{j=1}^{n-1} (1-\abs{\alpha_j}^2)
\end{equation}
But since $z\Phi_n = z^{n+1} + \text{ lower order}$, 
\begin{equation} \lb{4.17x}
\langle 1, z\Phi_n\rangle =\bar c_{n+1} + \text{ polynomial in } (c_0, c_1, \dots, 
c_n, \bar c_1, \dots, \bar c_n)
\end{equation}
This plus induction implies \eqref{4.15}. 

Given \eqref{4.14} and \eqref{4.15} plus the theorem of Schur that any 
$\{\gamma_j\}_{j=0}^{n-1}$ in $\bbD^n$ is allowed, and the theorem of Verblunsky 
that any $\{\alpha_j\}_{j=0}^{n-1}$ in $\bbD^n$ is allowed, we get 
\eqref{4.11A} inductively.

\subsection{Improved Exponential Decay Estimates} \lb{s4.5}

In \cite{NT89}, Nevai-Totik proved that 
\begin{equation} \lb{4.17}
\limsup_{n\to\infty}\, \abs{\alpha_n}^{1/n} = A<1 \Leftrightarrow d\mu_\s = 
0  \text{ and } D^{-1}(z)\text{ is analytic in } \{z\mid\abs{z}< A^{-1}\} 
\end{equation} 
providing a formula for the exact rate of exponential decay in terms of properties 
of $D^{-1}$. By analyzing their proof carefully, \cite[Section~7.2]{OPUC1} refines 
this to prove 

\begin{theorem}\lb{T4.8} Suppose 
\begin{equation} \lb{4.18}
\lim_{n\to\infty} \, \abs{\alpha_n}^{1/n} =A<1 
\end{equation}
and define 
\begin{equation} \lb{4.19}
S(z)=\sum_{n=0}^\infty \alpha_n z^n
\end{equation} 
Then $S(z) +\ol{D(1/\bar z)}\, D(z)^{-1}$ has an analytic continuation to 
$\{z\mid A<\abs{z} <A^{-2}\}$. 
\end{theorem} 

The point of this theorem is that both $S(z)$ and $\ol{D(1/\bar z)}\, D(z)^{-1}$  
have singularities on the circle of radius $A^{-1}$ ($S$ by \eqref{4.18} and 
$D^{-1}$ by \eqref{4.17}), so the fact that the combination has the continuation is 
a strong statement.  

Theorem~\ref{T4.8} comes from the same formula that Nevai-Totik \cite{NT89} use, 
namely, if $d\mu_\s =0$ and $\kappa_\infty =\prod_{n=0}^\infty (1-\abs{\alpha_n}^2)^{1/2}$, 
then 
\begin{equation} \lb{4.18x}
\alpha_n =-\kappa_\infty \int \ol{\Phi_{n+1}(e^{i\theta})}\, D(e^{i\theta})^{-1} \, 
d\mu(\theta)
\end{equation} 
We combine this with an estimate of Geronimus \cite{GBk1} that 
\begin{equation} \lb{4.19x}
\|\varphi_{n+1}^* -D^{-1}\|_{L^2 (\partial\bbD, d\mu)} \leq \sqrt2 
\biggl(\, \sum_{j=n+1}^\infty\, \abs{\alpha_j}^2\biggr)^{1/2}
\end{equation}
and $D^{-1}\, d\mu =\bar D\f{d\theta}{2\pi}$ to get 
\begin{equation} \lb{4.20}
\alpha_n +\sum_{j=n}^\infty d_{j,-1} \bar d_{j-n,1} = O((A^{-1} -\veps)^{-2n}) 
\end{equation}
where $D(z)=\sum_{j=0}^\infty d_{j,1}z^j$, $D(z)^{-1}=\sum_{j=0}^\infty d_{j,-1}z^j$. 
\eqref{4.20} is equivalent to analyticity of $S(z)+\ol{D(1/\bar z)}\, D(z)^{-1}$ in 
the stated region. 

One consequence of Theorem~\ref{T4.8} is  

\begin{corollary}\lb{C4.9} Let $b\in\bbD$. Then 
\begin{equation} \lb{4.21}
\f{\alpha_{n+1}}{\alpha_n} = b+O(\delta^n)
\end{equation}
for some $\delta<1$ if and only if $D^{-1}(z)$ is meromorphic in $\{z\mid \abs{z} < 
\abs{b}^{-1} + \delta'\}$ for some $\delta'$ and $D(z)^{-1}$ has only a single pole at 
$z=1/b$ in this disk. 
\end{corollary} 

This result is not new; it is proven by other means in Barrios-L\'opez-Saff \cite{BLS}. 
Our approach leads to a refined form of \eqref{4.21}, namely, 
\begin{equation} \lb{4.22}
\alpha_n =-Cb^n + O((b\delta)^n)
\end{equation} 
with 
\begin{equation} \lb{4.23}
C=\bigl[\, \lim_{z\to b^{-1}}\, (1-zb) D(z)^{-1} \bigr]\, \ol{D(\bar b)}
\end{equation}

One can get more. If $D(z)^{-1}$ is meromorphic in $\{z\mid \abs{z} 
<A^{-2}\}$, one gets an asymptotic expansion of $\alpha_n$ of the form 
\[
\alpha_n =\sum_{j=1}^\ell P_{m_j}(n) z_j^n +O((A^{-2}-\veps)^{-n}) 
\]
where the $z_j$ are the poles of $D^{-1}$ in $\{z\mid\abs{z}<A^{-2}\}$ and $P_{m_j}$ 
are polynomials of degree $m_j =$ the order of the pole at $m_j$. There are also results 
relating asymptotics of $\alpha_n$ of the form $\alpha_n =Cb^n n^k (1+o(1))$ to 
asymptotics of $d_{n,-1}$ or the form $d_{n,-1}=C_1 b^n n^k (1+o(1))$.

\subsection{Rakhmanov's Theorem on an Arc with Eigenvalues in the Gap} \lb{s4.6}

Rakhmanov \cite{Rakh77} proved a theorem that if \eqref{1.1} holds with $w(\theta)
\neq 0$ for a.e.~$\theta$, then $\lim_{n\to\infty} \abs{\alpha_n}=0$ (see also 
\cite{Rakh83,MNT85a,Nev91,Kh2000}). In \cite[Section~13.4]{OPUC2}, we prove the following 
new result related to this. Define for $a\in (0,1)$ and $\lambda\in\partial\bbD$ 
\begin{equation} \lb{4.23A}
\Gamma_{a,\lambda} =\{z\in\partial\bbD\mid\arg (\lambda z)>2\arcsin (a)\}
\end{equation}
and $\ess\,\supp(d\mu)$ of a measure as points $z_0$ with $\{z\mid\abs{z-z_0}<\veps\} 
\cap\supp(d\mu)$ an infinite set for all $\veps >0$. Then 

\begin{theorem}\lb{T4.10} Let $d\mu$ be given by \eqref{1.1} so that $\ess\,\supp(d\mu) 
=\Gamma_{a,\lambda}$ and $w(\theta) >0$ for a.e.~$e^{i\theta}\in\Gamma_{a,\lambda}$. 
Then 
\begin{equation} \lb{4.24}
\lim_{n\to\infty}\, \abs{\alpha_n(d\mu)} =a \qquad \lim_{n\to\infty} \,  
\ol{\alpha_{n+1}(d\mu)}\, \alpha_n (d\mu) = a^2\lambda
\end{equation}
\end{theorem} 

We note that, by rotation invariance, one need only look at $\lambda =1$. 
$\Gamma_a \cup\{1\}$ is known (Geronimus \cite{Ger66,Ger77}; see also 
\cite{Gol99,Gol2,GNA2,GNA1,PS1,PS2,PS3,PS4}) to be exactly the spectrum for 
$\alpha_n\equiv a$ and the spectrum on $\Gamma_a$ is purely a.c.~with 
$w(\theta)>0$ on $\Gamma_a^{\intt}$. 

Theorem~\ref{T4.10} can be viewed as a combination of two previous extensions of 
Rakhmanov's theorem. First, Bello-L\'opez \cite{BRLL} proved \eqref{4.24} if 
$\ess\,\supp (d\mu)=\Gamma_{a,\lambda}$ is replaced by $\supp(d\mu)=
\Gamma_{a,\lambda}$. Second, Denisov \cite{DenPAMS} proved an analog of Rakhmanov's 
theorem for OPRL. By the mapping of measures on $\partial\bbD$ to measures on 
$[-2,2]$ due to Szeg\H{o} \cite{Sz22a} and the mapping of Jacobi coefficients to 
Verblunsky coefficients due to Geronimus \cite{Ger46}, Rakhmanov's theorem immediately 
implies that if a Jacobi matrix has $\supp (d\gamma)=[-2,2]$ and $d\gamma=f(E)+ 
d\gamma_\s$ with $f(E)>0$ on $[-2,2]$, then $a_n\to 1$ and $b_n\to 0$. What Denisov 
\cite{DenPAMS} did is extend this result to only require $\ess\,\supp(d\gamma)=[-2,2]$. 

Thus, Theorem~\ref{T4.10} is essentially a synthesis of the Bello-L\'opez \cite{BRLL} 
and Denisov \cite{DenPAMS} results. One difficulty in such a synthesis is that Denisov 
relies on Sturm oscillation theorems and such a theorem does not seem to be applicable for 
OPUC. Fortunately, Nevai-Totik \cite{NTppt} have provided an alternate approach to 
Denisov's result using variational principles, and their approach --- albeit with some 
extra complications --- allows one to prove Theorem~\ref{T4.10}. The details are in 
\cite[Section~13.4]{OPUC2}.

\subsection{A Birman-Schwinger Principle for OPUC} \lb{s4.7}

Almost all quantitative results on the number of discrete eigenvalues for 
Schr\"odinger operators and OPRL depend on a counting principle of Birman 
\cite{Bir61} and Schwinger \cite{Schw}. In \cite[Section~10.15]{OPUC2}, we have 
found an analog for OPUC by using a Cayley transform and applying the 
Birman-Schwinger idea to it. Because of the need to use a point in $\partial\bbD$  
about which to base the Cayley transform, the constants that arise are not universal. 
Still, the method allows the proof of perturbation results like the following from 
\cite[Section~12.2]{OPUC2}: 

\begin{theorem}\lb{T4.11} Suppose $d\mu$ has Verblunsky coefficients $\alpha_j$ 
and there exists $\beta_j$ with $\beta_{j+p}=\beta_j$ for some $p$ and 
\begin{equation} \lb{4.25}
\sum_{j=0}^\infty j \abs{\alpha_j-\beta_j} <\infty 
\end{equation}
Then $d\mu$ has an essential support whose complement has at most $p$ gaps, and 
each gap has only finitely many mass points. 
\end{theorem} 

\begin{theorem}\lb{T4.12} Suppose $\alpha$ and $\beta$ are as in Theorem~\ref{T4.11}, 
but \eqref{4.25} is replaced by 
\begin{equation} \lb{4.26}
\sum_{j=0}^\infty\, \abs{\alpha_j-\beta_j}^p <\infty
\end{equation} 
for some $p\geq 1$. Then 
\begin{equation} \lb{4.27}
\sum_{z_j=\text{ mass points in gaps}}\, \dist (z_j, \ess\,\sup(d\mu))^q <\infty 
\end{equation}
where $q>\f12$ if $p=1$ and $q\geq p-\f12$ if $p>1$. 
\end{theorem} 

Theorem~\ref{T4.11} is a bound of Bargmann type \cite{Barg3}, while Theorem~\ref{T4.12} 
is of Lieb-Thirring type \cite{Lieb}. We have not succeeded in proving $q=\f12$ for $p=1$ 
whose analog is known for Schr\"odinger operators \cite{Weidl, HLT} and OPRL \cite{HunS}.

\subsection{Rotation Number for OPUC} \lb{s4.8}

Rotation numbers and their connection to the density of states have been an important 
tool in the theory of Schr\"odinger operators and OPRL (see Johnson-Moser \cite{JoMo}). 
Their analog for OPUC has a twist, as seen from the following theorem from 
\cite[Section~8.3]{OPUC1}: 

\begin{theorem}\lb{T4.13} $\arg (\Phi_n (e^{i\theta}))$ is monotone increasing 
in $\theta$ on $\partial\bbD$ and defines a measure $d\arg (\Phi_n (e^{i\theta}))/d\theta$ 
of total mass $2\pi n$. If the density of zeros $d\nu_n$ has a limit $d\nu$ supported 
on $\partial\bbD$, then 
\begin{equation} \lb{4.28}
\f{1}{2\pi n} \, \f{d\arg (\Phi_n (e^{i\theta}))}{d\theta}\to \f12\, d\nu 
+ \f12\, \f{d\theta}{2\pi}
\end{equation}
weakly. 
\end{theorem} 

Given the OPRL result, the $\f12 \, \f{d\theta}{2\pi}$ is surprising. In a sense, it comes 
from the fact that the transfer matrix obeys $\det (T_n)=z^n$ rather than determinant $1$. 
The proof of Theorem~\ref{T4.13} comes from an exact result that in turn comes from looking 
at $\arg (e^{i\theta}-z_0)$ for $z_0\in\bbD$: 
\begin{equation} \lb{4.29}
\f{1}{2\pi n} \, \f{d\arg (\Phi_n (e^{i\theta}))}{d\theta} = \f12\, \calP(d\nu_n) + 
\f12 \, \f{d\theta}{2\pi} 
\end{equation}
where $\calP$ is the dual of Poisson kernel viewed as a map of $C(\partial\bbD)$ to $C(\bbD)$, 
that is, 
\begin{equation} \lb{4.30}
\calP(d\gamma) = \f{1}{2\pi} \int \f{1-\abs{r}^2}{1+r^2 -2r\cos\theta}\, d\gamma 
(re^{i\theta})
\end{equation}

\bigskip
\section{Periodic Verblunsky Coefficients} \lb{s5} 

In this section, we describe some new results/approaches for Verblunsky coefficients 
$\{\alpha_n\}_{n=0}^\infty$ that obey 
\begin{equation} \lb{5.1} 
\alpha_{n+p}=\alpha_n 
\end{equation} 
for some $p$. We'll normally suppose $p$ is even. If it is not, one can use the fact that 
$(\alpha_0, 0, \alpha_1, 0, \alpha_2, 0,\dots)$ is the Verblunsky coefficients of the 
measure $\f12 d\mu (e^{2i\theta})$ and it has \eqref{5.1} with $p$ even, so one can 
read off results for $p$ odd from $p$ even. 

The literature is vast for Schr\"odinger operators with periodic potential called 
Hill's equation after Hill \cite{Hill1886}. The theory up to the 1950's is summarized 
in Magnus-Winkler \cite{MagWin} whose key tool is the discriminant; see also Reed-Simon 
\cite{RS4}. There was an explosion of ideas following the KdV revolution, including 
spectrally invariant flows and abelian functions on hyperelliptic Riemann surfaces. 
Key papers include McKean-van Moerbeke \cite{McvM1}, Dubrovin et al.~\cite{DubMatNov}, 
and Trubowitz \cite{Trub1}. Their ideas have been discussed for OPRL; see especially 
Toda \cite{Toda}, van Moerbeke \cite{vMoer}, and Flaschka-McLaughlin \cite{FlMcL}. 

For OPUC, the study of measures associated with \eqref{5.1} goes back to Geronimus 
\cite{Ger44} with a fundamental series of papers by Peherstorfer-Steinbauer 
\cite{Pe01,Pe03,PS1,PS2,PS3,PS4,PS6,PS5} and considerable literature on the case 
$p =1$ (i.e., constant $\alpha$); see, for example, 
\cite{Ger66,Ger77,Gol99,Gol2,GNA1,GNA2,Kh2000,Khr}. The aforementioned literature 
on OPUC used little from the the work on Hill's equation; work that does make a 
partial link is Geronimo-Johnson \cite{GJo1}, which discussed almost periodic 
Verblunsky coefficients using abelian functions. Simultaneous with our work 
reported here, Geronimo-Gesztesy-Holden \cite{GGH04} have discussed this further, 
including work on isospectral flows. Besides the work reported here, Nenciu-Simon 
\cite{NenSim} have found a symplectic structure on $\bbD^p$ for which the coefficients 
of the discriminant Poisson commute (this is discussed in \cite[Section~11.11]{OPUC2}.

\subsection{Discriminant and Floquet Theory} \lb{s5.1}

For Schr\"odinger operators, it is known that the discriminant is just the trace of the 
transfer matrix. Since the transfer matrix has determinant one in this case, the 
eigenvalues obey $x^2 -\tr(T)x+1=0$, which is the starting point for Floquet theory. 
For OPUC, the transfer matrix, $T_p(z)$, of \eqref{4.4a} has $\det (T_p(z))=z^p$, so 
it is natural to define the discriminant by 
\begin{equation} \lb{5.2} 
\Delta (z)=z^{-p/2} \tr (T_p(z)) 
\end{equation}
which explains why we take $p$ even. Because for $z=e^{i\theta}$, $A(\alpha,z)\in\bbU(1,1)$ 
(see \cite[Section~10.4]{OPUC2} for a discussion of $\bbU(1,1)$), $\Delta(z)$ is real on 
$\partial\bbD$, so 
\begin{equation} \lb{5.3} 
\Delta (1/\bar z)=\ol{\Delta(z)}
\end{equation}

Here are the basic properties of $\Delta$: 

\begin{theorem}\lb{T5.1} 
\begin{SL} 
\item[{\rm{(a)}}] All solutions of $\Delta(z)-w=0$ with $w\in (-2,2)$ are simple zeros 
and lie in $\partial\bbD$ {\rm{(}}so are $p$ in number{\rm{)}}. 

\item[{\rm{(b)}}] $\{z\mid \Delta(z) \in (-2,2)\}$ is $p$ disjoint intervals on $\partial\bbD$ 
whose closures $B_1, \dots, B_p$ can overlap at most in single points. The complements 
where $\abs{\Delta(z)}>2$ and $z\in\partial\bbD$ are ``gaps," at most $p$ in number.  

\item[{\rm{(c)}}] On $\cup B_j$, $d\mu$ is purely a.c.~{\rm{(}}i.e., in terms 
of \eqref{1.1}, $\mu_\s (\cup_{j=1}^p B_j)=0$ and $w(\theta)>0$ for a.e.~$\theta\in 
\cup_{j=1}^p B_j${\rm{)}}. 

\item[{\rm{(d)}}] $\mu\restriction (\partial\bbD\backslash \cup_{j=1}^p B_j)$ consists of 
pure points only with at most one pure point per gap. 

\item[{\rm{(e)}}] For all $z\in\bbC\backslash\{0\}$, the Lyapunov exponent $\lim_{n\to\infty} 
\|T_n(z)\|^{1/n}$ exists and obeys 
\begin{equation} \lb{5.4} 
\gamma(z)= \f12\, \log(z) + \f{1}{p}\, \log \biggl| \f{\Delta(z)}{2} + \sqrt{\f{\Delta^2}{4}-1}\,\biggr| 
\end{equation}
where the branch of square root is taken that maximizes the log. 

\item[{\rm{(f)}}] If $B=\cup_{j=1}^p B_j$, then the logarithmic capacity of $B$ is given by 
\begin{equation} \lb{5.5x} 
C_B = \prod_{j=0}^{p-1} (1-\abs{\alpha_j}^2)^{1/p}  
\end{equation}
and $-[\gamma(z) + \log C_B]$ is the equilibrium potential for $B$. 

\item[{\rm{(g)}}] The density of zeros is the equilibrium measure for $B$ and given in 
terms of $\Delta$ by 
\begin{equation} \lb{5.6x} 
d\nu(\theta) = V(\theta)\, \f{d\theta}{2\pi} 
\end{equation}
where $V(\theta)=0$ on $\partial\bbD\backslash \cup_{j=1}^p B_j$, and on $B_j$ is given by 
\begin{equation} \lb{5.7x} 
V(\theta)=\f{1}{p}\,\, \f{\abs{\Delta'(e^{i\theta})}}{\sqrt{4-\Delta^2 (e^{i\theta})}\,} 
\end{equation} 
where $\Delta' (e^{i\theta}) = \f{\partial}{\partial\theta}\, \Delta (e^{i\theta})$. 

\item[{\rm{(h)}}] $\nu (B_j) =1/p$ 
\end{SL} 
\end{theorem} 

For proofs, see \cite[Section~11.1]{OPUC2}. The proofs are similar to those for Schr\"odinger 
operators. That the density of zeros is an equilibrium measure has been emphasized by 
Saff, Stahl, and Totik \cite{SaffTot,SST}. While not expressed as the trace of a 
transfer matrix, $\Delta$ is related to the (monic) Tchebychev polynomial, $T$, of 
Peherstorfer-Steinbauer \cite{PS2} by 
\[
\Delta (z) =z^{-p/2} C_B^{-1/2} T(z) 
\]
and some of the results in Theorem~\ref{T5.1} are in their papers. 

One can also relate $\Delta$ to periodized CMV matrices, an OPUC version of Floquet theory. 
As discussed in Section~\ref{s3.7}, $\calE_p(\beta)$ is defined by restricting $\calE$ 
to sequences obeying $u_{n+p}=\beta u_n$ for all $n$. $\calE_p$ can be written as a $p\times p$ 
matrix with an $\calL\calM$ factorization. With $\Theta$ given by \eqref{3.4}, $\calE_p 
(\beta)=\calL_p\calM_p(\beta)$ 
\begin{equation} \lb{5.5}
\begin{aligned} 
\calM_p (\beta) &= \begin{pmatrix} 
-\alpha_{p-1} & {} & {} & {} & \rho_p \beta_{-1} \\ 
{} & \Theta_1 & {} & {} & {} \\
{} & {} & \ddots & {} & {}   \\
{} & {} & {} & \Theta_{p-3} & {} \\
\rho_{p-1}\beta & {} & {} & {} & \bar\alpha_{p-1} 
\end{pmatrix} \\ 
\\
\calL_p &= \begin{pmatrix} 
\Theta_0 & {} & {} & {} & {} \\
{} & \ddots & {} & {} & {} \\
{} & {} & \ddots & {} & {} \\
{} & {} & {} & {} & \ddots & {} \\
{} & {} & {} & {} & {} &\Theta_{p-2} 
\end{pmatrix}
\end{aligned} 
\end{equation}
Then: 

\begin{theorem}\lb{T5.2} 
\begin{SL} 
\item[{\rm{(a)}}] The following holds:
\begin{equation} \lb{5.6} 
\det(z-\calE_p(\beta)) = \prod_{j=0}^{p-1} (1-\abs{\alpha_j}^2)^{1/2p} z^{p/2} 
[\Delta(z)-\beta-\beta^{-1}] 
\end{equation}
\item[{\rm{(b)}}] $\calE$ is a direct integral of $\calE_p (\beta)$. 
\end{SL} 
\end{theorem}

\subsection{Generic Potentials} \lb{s5.2}

The following seems to be new; it is an analog of a result \cite{Sim75} for 
Schr\"odinger operators. 

Notice that for any $\{\alpha_j\}_{j=0}^p \in\bbD^p$, one can define a discriminant 
$\Delta(z,\{\alpha_j\}_{j=0}^{p-1})$ for the period $p$ Verblunsky coefficients that 
agree with $\{\alpha_j\}_{j=0}^{p-1}$ for $j=0,\dots, p-1$. 

\begin{theorem}\lb{T5.3} The set of $\{\alpha_j\}_{j=0}^p \in\bbD^p$ for which 
$\Delta(z)$ has all gaps open is a dense open set. 
\end{theorem} 

\cite{OPUC2} has two proofs of this theorem: one in Section~11.10 uses Sard's 
theorem and one is perturbation theoretic calculation that if 
$\abs{\tr(z,\{\alpha_j^{(0)}\}_{j=0}^{p-1})}=2$, then 
$\abs{\tr(z,\{\alpha_j^{(0)} + (e^{i\eta}-1) \delta_{jk} \alpha_k^{(0)}\}_{j=0}^{p-1})} 
=2+2\eta^2 (\rho_k^{(0)})^2 \abs{\alpha_k^{(0)}} + O(\eta^3)$.

\subsection{Borg's Theorems} \lb{s5.3}

In \cite{Borg1}, Borg proved several theorems about the implication of closed gaps. 
Further developments of Borg's results for Schr\"odinger equations or for OPRL are 
in Hochstadt \cite{Hoch1,Hoch2,Hoch,Hoch3,Hoch4,Hoch5,Hoch6,Hoch7}, Clark et al.~\cite{CGHL}, 
Trubowitz \cite{Trub1}, and Flaschka \cite{Fla3}. In \cite[Section~11.14]{OPUC2}, we 
prove the following analogs of these results: 

\begin{theorem} \lb{T5.4} If $\{\alpha_j\}_{j=0}^\infty$ is a periodic sequence of Verblunsky 
coefficients so $\supp(d\mu)=\partial\bbD$ {\rm{(}}i.e., all gaps are closed{\rm{)}}, then 
$\alpha_j\equiv 0$. 
\end{theorem}

\cite{OPUC2} has three proofs of this: one uses an analog of a theorem of Deift-Simon 
\cite{S169} that $d\mu/d\theta \geq 1/2\pi$ on the essential support of the a.c.~spectrum 
of any ergodic system, one tracks zeros of the Wall polynomials, and one uses the analog 
of Tchebychev's theorem for the circle that any monic Laurent polynomial real on $\partial 
\bbD$ has $\max_{z\in\partial\bbD}\abs{L(z)}\geq 2$. 

\begin{theorem}\lb{T5.5} If $p$ is even and $\{\alpha_j\}_{j=0}^\infty$ has period  
$2p$, then if all gaps with $\Delta(z)=-2$ are closed, we have $\alpha_{j+p}=\alpha_j$,  
and if all gaps with $\Delta(z)=2$ are closed, then $\alpha_{j+p}=-\alpha_j$. 
\end{theorem} 

\begin{theorem}\lb{T5.6} Let $p$ be even and suppose for some $k$ that $\alpha_{kp+j} 
=\alpha_j$ for all $j$. Suppose for some labelling of $\{w_j\}_{j=0}^{k_p-1}$ of the zeros 
of the derivative $\partial\Delta/\partial\theta$ labelled counterclockwise, we have 
$\abs{\Delta (w_j)}=2$ if $j\not\equiv 0 \text{\rm{ mod }}k$. Then $\alpha_{p+j}= 
\omega\alpha_j$ where $\omega$ is a $k$-th root of unity. 
\end{theorem} 

The proof of these last two theorems depends on the study of the Carath\'eodory function 
for periodic Verblunsky coefficients as meromorphic functions on a suitable hyperelliptic 
Riemann surface.

\subsection{Green's Function Bounds} \lb{s5.4}

In \cite[Section~10.14]{OPUC2}, we develop the analog of the Combes-Thomas \cite{ComTh} 
method for OPUC and prove, for points in $\partial\bbD\backslash\supp(d\mu)$, the 
Green's function (resolvent matrix elements of $(\calC-z)^{-1}$ with $\calC$ the 
CMV matrix) decays exponentially in $\abs{n-m}$. The rate of decay in these estimates 
goes to zero at a rate faster than expected in nice cases. For periodic Verblunsky 
coefficients, one expects behavior similar to the free case for OPRL or Schr\"odinger 
operators --- and that is what we discuss here. An energy $z_0\in\partial\bbD$ at 
the edge of a band is called a resonance if $\sup_n \abs{\varphi_n (z_0)}<\infty$. 
For the family of measures, $d\mu_\lambda$, with Verblunsky coefficients $\alpha_n = 
\lambda\alpha_n^{(0)}$ and a given $z_0$, there is exactly one $\lambda$ for which 
$z_0$ is a resonance (for the other values, $\varphi_n (z_0,d\mu_\lambda)$ grows 
linearly in $n$). Here is the bound we prove in \cite[Section~11.12]{OPUC2}: 

\begin{theorem}\lb{T5.7} Let $\{\alpha_n\}_{n=0}^\infty$ be a periodic family of 
Verblunsky coefficients. Suppose $G=\{z=e^{i\theta}\mid \theta_0 < \theta < \theta_1\}$ 
is an open gap and $e^{i\theta_0}$ is not a resonance. Let 
\[
G_{nm}(z) = \langle\delta_n, (\calC(\alpha) - z)^{-1}\delta_m\rangle
\]
Then for $z=e^{i\theta}$ with $z\in G$ and $\abs{\theta - \theta_0} < 
\abs{\theta-\theta_1}$, we have
\begin{align*}  
\sup_{n,m}\, \abs{G_{nm}(z)} &\leq C_1 \abs{z-e^{i\theta_0}}^{-1/2}  \\ 
\sup_{\text{\rm{such} } z}\, \abs{G_{nm}(z)} &\leq C_2 (n+1)^{1/2} (m+1)^{1/2} 
\end{align*} 
and similarly for $z$ approaching $e^{i\theta_1}$. 
\end{theorem}  

The proof depends on bounds on polynomials in the bands of some independent interest. 

\begin{theorem}\lb{T5.8}  Let $\{\alpha_n\}_{n=0}^\infty$ be a sequence of 
periodic Verblunsky coefficients, and let $B^\intt$ be the union of the interior 
of the bands. Let $\calE_1$ be the set of band edges by open gaps and $\calE_2$ 
the set of band edges by closed gaps. Define 
\[
d(z) =\min (\dist(z,\calE_1), \dist (z,\calE_2)^2) 
\]
Then 
\begin{alignat*}{2} 
&(1) \qquad  \sup_n \, \abs{\varphi_n(z)} &&\leq C_1\, d(z)^{-1/2}  \\ 
&(2) \qquad \sup_{z\in B^\intt}\, \abs{\varphi_n(z)} && \leq C_2 n 
\end{alignat*} 
where $C_1$ and $C_2$ are $\{\alpha_n\}_{n=0}^\infty$ dependent constants. 
\end{theorem} 

{\it Remark.}  One can, with an extra argument, show $d(z)$ can be replaced 
by $\dist (z,\calE_1)$ which differs from $d(z)$ only when there is a closed gap. 
That is, there is no singularity in $\sup_n \abs{\varphi_n(z)}$ at band edges 
next to closed gaps.

\subsection{Isospectral Results} \lb{s5.5}

In \cite[Chapter~11]{OPUC2}, we prove the following theorem: 

\begin{theorem}\lb{T5.9} Let $\{\alpha_j\}_{j=0}^{p-1}$ be a sequence in $\bbD^p$ so 
$\Delta(z, \{\alpha_j\}_{j=0}^{p-1})$ has $k$ open gaps. Then 
$\{\{\beta_j\}_{j=0}^{p-1}\in\bbD^p\mid\Delta(z,\{\beta_j\}_{j=0}^{p-1}) =
\Delta(z,\{\alpha_j\}_{j=0}^{p-1})\}$ is a $k$-dimensional torus. 
\end{theorem}

This result for OPUC seems to be new, although its analog for finite-gap Jacobi matrices 
and Schr\"odinger operators (see, e.g., \cite{McvM1,DubMatNov, vMoer}) is well known 
and it is related to results on almost periodic OPUC by Geronimo-Johnson \cite{GJo1}. 

There is one important difference between OPUC and the Jacobi/Schr\"odinger case. In the  
later, the infinite gap doesn't count in the calculation of dimension of torus, so the 
torus has a dimension equal to the genus of the Riemann surface for the $m$-function. 
In the OPUC case, all gaps count and the torus has dimension one more than the genus. 

The torus can be defined explicitly in terms of natural additional data associated to 
$\{\alpha_j\}_{j=0}^{p-1}$. One way to define the data is to analytically continue the 
Carath\'eodory function, $F$, for the periodic sequence. One cuts $\bbC$ on the 
``combined bands," that is, connected components of $\{e^{i\theta}\mid\abs{\Delta (e^{i\theta})} 
\leq 2\}$, and forms the two-sheeted Riemann surface associated to $\sqrt{\Delta^2 -4}$. 
On this surface, $F$ is meromorphic with exactly one pole on each ``extended gap." By 
extended gap, we mean the closure of the two images of a gap on each of two sheets 
of the Riemann surface. The ends of the gap are branch points and join the two images 
into a circle. The $p$ points, one on each gap, are thus $p$-dimensional torus, and 
the refined version of Theorem~\ref{T5.9} is that there is exactly one Carath\'eodory 
function associated to a period $p$ set of Verblunsky coefficients with specified poles. 

Alternately, the points in the gaps are solutions of $\Phi_p(z)-\Phi_p^*(z)=0$ with 
sheets determined by whether the points are pure points of the associated measure or not. 

\cite{OPUC2} has two proofs of Theorem~\ref{T5.9}: one using the Abel map on the 
above referenced Riemann surface and one using Sard's theorem.

\subsection{Perturbation Conjectures}\lb{s5.6}

\cite{OPUC1,OPUC2} have numerous conjectures and open problems. We want to end this 
section with a discussion of conjectures that describe perturbations of periodic Verblunsky 
coefficients. We discuss the Weyl-type conjecture in detail. As a model, consider 
Theorem~\ref{T3.4} when $\alpha_n\equiv a\neq 0$. For $\ess\,\supp(d\nu)$ to be 
$\Gamma_{a,1}$, the essential support for $\alpha_n\equiv a$, it suffices that 
$\abs{\alpha_n}\to a$ and $\alpha_{n+1}/\alpha_n\to 1$. This suggests 

\begin{conjecture} \lb{Con5.10} Fix a period $p$ set of Verblunsky coefficients with 
discriminant $\Delta$. Let $M$ be the set of period $p$ (semi-infinite) sequences 
with discriminant $\Delta$ and let $S\subset\partial\bbD$ be their common 
essential support. Suppose 
\[
\lim_{j\to\infty}\, \inf_{\alpha\in  M}\, \biggl[ \, \sum_{n=1}^\infty 
e^{-n} \abs{\beta_{j+n}-\alpha_n}\biggr]=0
\]
Then if $\nu$ is the measure with Verblunsky coefficients $\beta$, then $\ess\,\supp 
(d\nu)=S$. 
\end{conjecture} 

Thus, limit results only hold in the sense of approach to the isospectral manifold. 
There are also conjectures in \cite{OPUC2} for extensions of Szeg\H{o}'s and 
Rakmanov's theorems in this context.

\bigskip
\section{Spectral Theory Examples} \lb{s6} 

\cite[Chapter 12]{OPUC2} is devoted to analysis of specific classes of Verblunsky 
coefficients, mainly finding analogs of known results for Schr\"odinger or 
discrete Schr\"odinger equations. Most of these are reasonably straightforward, 
but there are often some extra tricks needed and the results are of interest.

\subsection{Sparse and Decaying Random Verblunsky Coefficients} \lb{s6.1}

In \cite{S265}, Kiselev, Last, and Simon presented a thorough analysis of continuum 
and discrete Schr\"odinger operators with sparse or decaying random potentials, 
subjects with earlier work by Pearson \cite{Pear78}, Simon \cite{Sim157}, Delyon 
\cite{Dely,Sim178}, and Kotani-Ushiroya \cite{KU}. In \cite[Sections~12.3 and 12.7]{OPUC2}, 
I have found analogs of these results for OPUC: 

\begin{theorem}\lb{T6.1} Let $d\mu$ have the form \eqref{1.1}. Let $\{n_\ell\}_{\ell=1}^\infty$ 
be a monotone sequence of positive integers with $\liminf_{\ell\to\infty} \f{n_{\ell+1}}{n_\ell} 
>1$ and 
\begin{equation} \lb{6.1} 
\alpha_j (d\mu) =0 \qquad\text{if } j\notin \{n_\ell\} 
\end{equation}
and 
\begin{equation} \lb{6.2} 
\sum_{j=0}^\infty \, \abs{\alpha_j (d\mu)}^2 <\infty 
\end{equation}
Then $\mu_\s =0$, $\supp(d\mu)=\partial\bbD$, and $w,w^{-1}\in\cap_{p=1}^\infty 
L^p (\partial\bbD, \f{d\theta}{2\pi})$.  
\end{theorem} 

This result was recently independently obtained by Golinskii \cite{Gol2003}: 

\begin{theorem}\lb{T6.2} Let $\{n_\ell\}_{\ell=1}^\infty$ be a monotone sequence of 
positive integers with $\lim \f{n_{\ell+1}}{n_\ell} =\infty$ so that \eqref{6.1} 
holds. Suppose $\lim_{j\to\infty} \abs{\alpha_j(d\mu)}=0$ and \eqref{6.2} fails. 
Then $d\mu$ is purely singular continuous. 
\end{theorem} 

\begin{theorem}\lb{T6.3} Let $\{\alpha_j(\omega)\}_{j=0}^\infty$ be a family of 
independent random variables with values in $\bbD$ with 
\begin{equation} \lb{6.3} 
\bbE (\alpha_j (\omega)) =0 
\end{equation}
and 
\begin{equation} \lb{6.4} 
\sum_{j=0}^\infty \bbE (\abs{\alpha_j(\omega)}^2) <\infty 
\end{equation}
Let $d\mu_\omega$ be the measure with $\alpha_j (d\mu_\omega) =\alpha_j(\omega)$. 
Then for a.e.~$\omega$, $d\mu_\omega$ has the form \eqref{1.1} with $d\mu_{\omega,\s}$ 
and $w(\theta)>0$ for a.e.~$\theta$.  
\end{theorem} 

This result is not new; it is a result of Teplyaev, with earlier results of Nikishin 
\cite{Nik87} (see Teplyaev \cite{Tep89,Tep91,Tep92,Tep94}). We state it for comparison 
with the next two theorems.

The theorems assume \eqref{6.3} and also 
\begin{alignat}{2}
& \sup_{\omega,j}\, \abs{\alpha_j(\omega)}<1 \quad \sup_\omega \, \abs{\alpha_j(\omega)} 
\to 0 \qquad &&\text{as } j\to\infty \lb{6.5} \\
& \bbE(\alpha_j(\omega)^2) =0 \lb{6.6} \\ 
& \bbE(\abs{\alpha_j(\omega)}^2)^{1/2} = \Gamma j^{-\gamma} \qquad &&\text{if } j>J_0 \lb{6.7} 
\end{alignat}

\begin{theorem}\lb{T6.4} If $\{\alpha_j(\omega)\}_{j=0}^\infty$ is a family of 
independent random variables so \eqref{6.3}, \eqref{6.5}, \eqref{6.6}, and \eqref{6.7} 
hold and $\Gamma >0$, $\gamma <\f12$, then for a.e.~pairs $\omega$ and $\lambda\in\partial 
\bbD$, $d\mu_{\lambda,\omega}$, the measure with $\alpha_j (d\mu_{\lambda,\omega}) = 
\lambda\alpha_j(\omega)$, is pure point with support equal to $\partial\bbD$ {\rm{(}}i.e., 
dense mass points{\rm{)}}. 
\end{theorem} 

\begin{theorem}\lb{T6.5} If $\{\alpha_j(\omega)\}_{j=0}^\infty$ is a family of 
independent random variables so \eqref{6.3}, \eqref{6.5}, \eqref{6.6}, and \eqref{6.7} 
hold for $\Gamma >0$, $\gamma =\f12$, and 
\begin{equation} \lb{6.8} 
\sup_{n,\omega}\, n^{1/2} \abs{\alpha_n (\omega)} <\infty  
\end{equation}
Then 
\begin{SL} 
\item[{\rm{(i)}}] If $\Gamma^2 >1$, then for a.e.~pairs $\lambda\in\partial\bbD$, 
$\omega\in\Omega$, $d\mu_{\lambda,\omega}$ has dense pure point spectrum. 
\item[{\rm{(ii)}}] If $\Gamma^2 \leq 1$, then for a.e.~pairs $\lambda\in\partial\bbD$, 
$\omega\in\Omega$, $d\mu_{\lambda,\omega}$ has purely singular continuous 
spectrum of exact Hausdorff dimension $1-\Gamma^2$ in that $d\mu_{\lambda,\omega}$ 
is supported on a set of dimension $1-\Gamma^2$ and gives zero weight to any set 
$S$ with $\dim (S) <1-\Gamma^2$. 
\end{SL} 
\end{theorem}  

For the last two theorems, a model to think of is to let $\{\beta_n\}_{n=0}^\infty$ 
be identically distributed random variables on $\{z\mid \abs{z}\leq r\}$ for some 
$r<1$ with a rotationally invariant distribution and to let $\alpha_n=\Gamma^{1/2} 
\bbE (\abs{\beta_1}^2)^{-1/2} \max (n,1)^{-\gamma}\beta_n$. 

The proofs of these results exploit Pr\"ufer variables, which for OPUC go back to 
Nikishin \cite{Nik85} and Nevai \cite{Nev87}.

\subsection{Fibonacci Subshifts} \lb{s6.2}

For discrete Schr\"odinger operators, there is an extensive literature 
\cite{Aubry80,KKT,OPRSS,Suto,BIST,BGH,DLCMP99,DLLMP99,DKL00,DL03,LZ02A,LZ02B} 
on subshifts (see \cite{OPUC2,Quef,Loth} for a definition of subshifts). In 
\cite[Section~12.8]{OPUC2}, we have analyzed the OPUC analog  of the most heavily 
studied of these subshifts, defined as follows: Pick $\alpha,\beta\in\bbD$. Let 
$F_1=\alpha$, $F_2=\alpha\beta$, and $F_{n+1}=F_nF_{n-1}$ for $n=2,3,\dots$. 
$F_{n+1}$ is a sequence which starts with $F_n$ and so there is a limit 
$F=\alpha,\beta,\alpha,\alpha,\beta,\alpha,\beta,\alpha,\alpha,\beta, \alpha, 
\alpha,\dots$. We write $F(\alpha,\beta)$ when we want to vary $\alpha$ and $\beta$. 

\begin{theorem}\lb{T6.6} The essential support of the measure $\mu$ with  
$\alpha_{\boldsymbol{\cdot}} (d\mu)=F(\alpha, \beta)$ is a closed perfect 
set of Lebesgue measure zero for any $\alpha\neq\beta$. For fixed $\alpha_0,\beta_0$ 
and a.e.~$\lambda\in\partial\bbD$, the measure with $\alpha_{\boldsymbol{\cdot}} 
(d\mu) =F(\lambda\alpha_0, \lambda\beta_0)$ is a pure point measure, with each 
pure point isolated and the limit points of the pure points a perfect set of 
$\f{d\theta}{2\pi}$-measure zero. 
\end{theorem} 

The proof follows that for Schr\"odinger operators with a few additional tricks needed.

\subsection{Dense Embedded Point Spectrum} \lb{s6.3}

Naboko \cite{Nab1,Nab2,Nab5,Nab3,Nab4} and Simon \cite{S254} constructed Schr\"odinger 
operators $-\f{d^2}{dx^2} +V(x)$ with $V(x)$ decaying only slightly slower than 
$\abs{x}^{-1}$ so there is dense embedded point spectrum. Naboko's method extends 
to OPUC.

\begin{theorem}\lb{T6.7} Let $g(n)$ be an arbitrary function with $0<g(n)\leq 
g(n+1)$ and $g(n)\to\infty$ as $n\to\infty$. Let $\{\omega_j\}_{j=0}^\infty$ be an 
arbitrary subset of $\partial\bbD$ which are multiplicatively rationally independent, 
that is, for no $n_1, n_2, \dots, n_k\in\bbZ$ other than $(0,0,\dots, 0)$, is it 
true that $\prod_{j=1}^k (\omega_j \omega_0^{-1})^{n_j}=1$. Then there exists a sequence 
$\{\alpha_j\}_{j=0}^\infty$ of Verblunsky coefficients with  
\[
\abs{\alpha_n} \leq \f{g(n)}{n}
\]
for all $n$ so that the measure $d\mu$ with $\alpha_j (d\mu) =\alpha_j$ has pure points at 
each $\omega_j$. 
\end{theorem} 

{\it Remark.} If $g(n)\leq n^{1/2-\veps}$, then $\abs{\alpha_n}\in\ell^2$ so, by 
Szeg\H{o}'s theorem, $d\mu$ has the form \eqref{1.1} with $w(\theta)>0$ for a.e.~$\theta$, 
that is, the point masses are embedded in a.c.~spectrum.

\subsection{High Barriers} \lb{s6.4} 

Jitomirskaya-Last \cite{JL1} analyzed sparse high barriers to get discrete 
Schr\"odinger operators with fractional-dimensional spectrum. Their methods can 
be applied to OPUC. Let $0<a<1$ and 
\begin{align} 
L &= 2^{n^n} \lb{6.9} \\ 
\alpha_j &= (1-\rho_j^2)^{-1/2} \lb{6.10} \\ 
\rho_j &= \begin{cases} L_n^{-(1-a)/2a} & j=L_n \\
0 & \text{otherwise} \end{cases} \lb{6.11} 
\end{align} 

\begin{theorem}\lb{T6.8} Let $\alpha_j$ be given by \eqref{6.10}/\eqref{6.11} and let 
$d\mu_\lambda$ be the Aleksandrov measures with $\alpha_j (d\mu_\lambda) = 
\lambda\alpha_j$. Then for Lebesgue a.e.~$\lambda$, $d\mu_\lambda$ has exact dimension 
$a$ in the sense that $d\mu_\lambda$ is supported on a set of Hausdorff dimension 
$a$ and gives zero weight to any set $B$ of Hausdorff dimension strictly less 
than $a$. 
\end{theorem}

\bigskip


\end{document}